\documentclass[twoside]{article}
\usepackage{graphicx,amssymb,mathrsfs,amsmath,color,fancyhdr,amsthm}
\usepackage[all]{xy}
\input cyracc.def

\renewcommand{\baselinestretch}{1.06} 
\catcode`@=11 \long\def\@makefntext#1{\noindent #1}
\newskip\tabcentering \tabcentering=1000pt plus 1000pt minus 1000pt
\def\REF#1{\par\hangindent\parindent\indent\llap{#1\enspace}}
\def\MCH#1#2{\setbox0=\hbox{\raise#1\hbox{#2}}\smash{\box0}}

\def\dl{\displaystyle}
\let\@oddfoot\@empty  \let\@evenfoot\@empty

\def\@evenhead{}\def\@oddhead{}
\def\@evenhead{\vbox{\hbox to \textwidth{\footnotesize\rm\hbox to
1.0cm{\thepage\hfill} \hfill\hspace{2mm}\footnotesize{
\emph{Xiao J Zhao M}}}}}
\def\@oddhead{\vbox{\hbox to \textwidth{\footnotesize
{\it A parameterization of the canonical bases} \hfill{\ } \hfill\hbox to
1cm{\hfill\thepage}}}}

\def\sec#1{\vspace{2mm}\noindent{{\bf #1}}\vspace{0.5mm}}
\def\subsec#1{\vspace{2mm}\leftline{\bf #1}} 
\def\th#1{\vspace{1mm}\noindent{\bf #1}\quad } 
\def\pf#1{\vspace{1mm}\noindent{\it #1}\quad}

\def\leq{\leqslant}
\def\geq{\geqslant}
       
  \def\hml{\end{document}}  \newsymbol\wjzhml 203F \def\no{\noindent}

\begin{document}
\abovedisplayskip=3pt plus 1pt minus 1pt 
\belowdisplayskip=3pt plus 1pt minus 1pt 

\def\le{\leqslant}
\def\ge{\geqslant}
\def\dl{\displaystyle}




\vspace{8true mm}

\renewcommand{\baselinestretch}{1.9}\baselineskip 19pt

\noindent{\LARGE\bf A parameterization of the canonical bases of affine modified quantized enveloping algebras}

\vspace{0.5 true cm}

\noindent{\normalsize\sf XIAO Jie$^{1}$ \& ZHAO
Minghui$^{1\,\dag}$ \footnotetext{\baselineskip 10pt
$^\dag$ Corresponding author\\
This work was supported by NSF of China (No. 11131001).}}

\vspace{0.2 true cm}
\renewcommand{\baselinestretch}{1.5}\baselineskip 12pt
\noindent{\footnotesize\rm $^1$Department of Mathematical Sciences, Tsinghua University, Beijing 100084, P. R. China \\
(email: jxiao@math.tsinghua.edu.cn, zhaomh08@mails.tsinghua.edu.cn)\vspace{4mm}}

\baselineskip 12pt \renewcommand{\baselinestretch}{1.18}
\noindent{{\bf Abstract}\small\hspace{2.8mm} 
For symmetrizable Kac-Moody Lie algebra $\textbf{g}$, Lusztig introduced the modified quantized enveloping algebra $\dot{\textbf{U}}(\textbf{g})$ and its canonical basis in [12]. In this paper, for finite and affine type symmetric Lie algebra $\textbf{g}$ we define a set which depend only on the root category and prove that there is a bijection between the set and the canonical basis of $\dot{\textbf{U}}(\textbf{g}),$ where the root category  is the $T^2$-orbit category of the derived category of Dynkin or tame quiver. Our method bases on one theorem of Lin, Xiao and Zhang in [9], which gave the  PBW-basis of $\textbf{U}^+(\textbf{g})$.

 }

\vspace{1mm} \no{\footnotesize{\bf Keywords:\hspace{2mm}
Ringel-Hall algebra, modified quantized enveloping algebra, canonical basis
}}

\no{\footnotesize{\bf MSC(2000):\hspace{2mm}16G20, 17B37 }
 \vspace{2mm}
\baselineskip 15pt
\renewcommand{\baselinestretch}{1.22}
\parindent=10.8pt  
\rm\normalsize\rm

\sec{1\quad Introduction}


Let $\textbf{U}^+$ be the positive part of the quantized enveloping algebra $\textbf{U}$ associated to a Cartan datum.
In the case of finite type, Lusztig give two approach to construct the canonical basis ([10]).
The first is an elementary algebraic construction. By using Ringel-Hall algebra realization to $\textbf{U}^+$,
the isomorphism classes of representations of the Dynkin quiver give a PBW-type basis of $\textbf{U}^+$ and there is an order on this basis.
Under this order, the transition matrix between this basis and a monomial basis is a triangular matrix with diagonal entries equal to $1$.
By a standard linear algebra method one can get a bar invariant basis, which is the canonical basis. The second is a geometric construction.
Lusztig construct the canonical basis by using perverse sheaves and intersection cohomology.
Then the geometric construction of the canonical basis was generalized to the cases of all type in [11].
In the case of affine type, Lin, Xiao and Zhang in [9] provide a process to construct a PBW-type basis and get the canonical basis by using Ringel-Hall algebra approach and the representations of tame quivers ([9]).


Let $\dot{\textbf{U}}$ be the modified quantized enveloping algebra obtained from $\textbf{U}$ by modifying the Cartan part $\textbf{U}^0$
to $\oplus_{\lambda\in P}\mathbb{Q}(v)\mathbf{1}_{\lambda}$, where $P$ is the weight lattice. This algebra has the same representations
with $\textbf{U}$. Lusztig consider it as the limit of the tensor product of the highest weight modules and lowest weight modules.
Then Lusztig define the canonical bases of the tensor products and the canonical basis on $\dot{\textbf{U}}$ can be obtained from them ([12][13]).
Kashiwara also study the algebra $\dot{\textbf{U}}$ and its canonical basis ([7]).

In [5][6], Happle study the derived category $D^b(\Lambda)$ of a finite dimensional algebra $\Lambda$, which is a triangulated category with the translation functor $T$. He found  that there is a bijection between the isomorphism classes of the indecomposable objects in $\mathcal{R}=D^b(\Lambda)/T^2$ and the root of the corresponding Lie algebra, if $\Lambda$ is hereditary and representation-finite, where $\mathcal{R}=D^b(\Lambda)/T^2$ is the $T^2$-orbit category of $D^b(\Lambda)$. It was proved in [14] that $\mathcal{R}=D^b(\Lambda)/T^2$ is still a triangulated category, so it is called a root category. In [14][15], Peng and Xiao recover the whole Lie algebra structure via the root categories of finite dimensional hereditary algebras.

Note that the canonical basis of $\dot{\textbf{U}}$ is defined abstractly and depend on the canonical basis of $\textbf{U}^+$. Inspired by the method of constructing Lie algebra from root category, we want to get a connection between the canonical basis of $\dot{\textbf{U}}$ and the objects of the root category using the PBW-type basis defined by Lin, Xiao and Zhang in the case of affine type. In this paper, we associate a set $\tilde{\mathcal{M}}$ to the root category. The definition of the set bases on the corresponding set defined by Lin, Xiao, and Zhang in [9] of a hereditary category. But the set depend only on the root category and not on the embedding of the hereditary category to the root category. Then for a fixed embedding of the hereditary category to the root category, we can get a bijection between the set $\tilde{\mathcal{M}}$ and the canonical basis of $\dot{\textbf{U}}\mathbf{1}_{\lambda}$ for every $\lambda\in P$. Hence we say that the set $\tilde{\mathcal{M}}$ we construct from the root category provides a parameterization of the canonical basis of $\dot{\textbf{U}}$.

Since [20], it has been an open problem: how to realize the whole quantized enveloping algebra by using Hall algebra from derived category or root category. A lot of effort has been paid on the progress ([2][8][19][21]) and the most recent progress is given by Bridgeland in [1]. We hope that the main result in the present paper can provide a strong evidence for the connection between the canonical basis and root category.

In Section 2, we first give the basic notations of quantized enveloping algebra and modified quantized enveloping algebra. Then we recall the definition of the Ringel-Hall algebra and root category. In Section 3, we consider the case of finite type, which is simpler and can reflect the idea clearly. In Section 4, we consider the case of affine case. We first recall the definition of the PBW-type basis of $\textbf{U}^+$. Then we define a set $\tilde{\mathcal{M}}$ on the root category and define a PBW-type basis of $\dot{\textbf{U}}\mathbf{1}_{\lambda}$ with index in $\tilde{\mathcal{M}}$. By the standard linear algebra method, we get a bar-invariant basis and prove that this is the leading term of the canonical basis. Hence we get a parameterization of the canonical basis by using the root category.

\sec{2\quad Preliminaries}

\subsec{2.1\quad Quantized enveloping algebra}

Let $\mathbb{Q}$ be the field of rational numbers and $\mathbb{Z}$ the ring of integers. Let $I$ be a finite index set with $|I|=n$ and $A=(a_{ij})_{i,j\in I}$ be a generalized Cartan matrix. Denote by $r(A)$ the rank of $A$. Let $P^{\vee}$ be a free abelian group of rank $2n-r(A)$ with a $\mathbb{Z}$-basis $\{h_i|i\in I\}\cup\{d_s|s=1,\ldots, n- r(A)\}$ and let $\mathfrak{h}=\mathbb{Q}\otimes_{\mathbb{Z}}P^{\vee}$ be the $\mathbb{Q}$-linear space spanned by $P^{\vee}$. We call $P^{\vee}$ the dual weight lattice and $\mathfrak{h}$ the Cartan subalgebra. We also define the weight lattice to be $P=\{\lambda\in\mathfrak{h}^{\ast}|\lambda(P^{\vee})\subset\mathbb{Z}\}$.

Set $\Pi^{\vee}=\{h_i|i\in I\}$ and choose a linearly independent subset $\Pi=\{\alpha_i|i\in I\}\subset\mathfrak{h}^{\ast}$ satisfying $\alpha_j(h_i)=a_{ij}$ and $\alpha_j(d_s)=0$ or $1$ for $i,j\in I$, $s=1,\ldots, n-\textrm{rank} A$. The elements of $\Pi$ are called simple roots, and the elements of $\Pi^{\vee}$ are called simple coroots. The quintuple $(A,\Pi,\Pi^{\vee},P,P^{\vee})$ is called a Cartan datum associated with the generalized Cartan matrix $A$.

We recall the definition of the quantized enveloping algebra. Assume that  $A=(a_{ij})_{i,j\in I}$ is a symmetric generalized Cartan matrix.

Fix an indeterminate $v$. For $n\in \mathbb{Z}$, we set
\begin{displaymath}
[n]_{v}=\frac{v^n-v^{-n}}{v-v^{-1}},
\end{displaymath}
and $[0]_{v}!=1$,  $[n]_{v}!=[n]_{v}[n-1]_{v}\cdots[1]_{v}$ for $n\in\mathbb{Z}_{>0}$. For nonnegative integers $m\geq n\geq 0$, the analogues of binomial coefficients are given by
\begin{displaymath}
\left[{m \atop n}\right]_v=\frac{[m]_v!}{[n]_v![m-n]_v!}.
\end{displaymath}
Then $[n]_v$ and $\left[{m \atop n}\right]_v$ are elements of the field $\mathbb{Q}(v)$.

The quantized enveloping algebra $\mathbf{U}$ associated with a Cartan datum $(A,\Pi,\Pi^{\vee},P,P^{\vee})$ is an associative algebra over $\mathbb{Q}(v)$ with $\mathbf{1}$ generated by the elements $E_i$, $F_i(i\in I)$ and $K_{\mu}(\mu\in P^{\vee})$ subject to the following relations:

(1) $K_{0}=\mathbf{1}$, $K_{\mu}K_{\mu'}=K_{\mu+\mu'}$ for all $\mu,\mu'\in P^{\vee}$;

(2) $K_{\mu}E_{i}K_{-\mu}=v^{\alpha_i(\mu)}E_i$ for all $i\in I$, $\mu\in P^{\vee}$;

(3) $K_{\mu}F_{i}K_{-\mu}=v^{-\alpha_i(\mu)}E_i$ for all $i\in I$, $\mu\in P^{\vee}$;

(4) $E_iF_i-F_iE_i=\delta_{ij}\frac{{K}_{i}-{K}_{-i}}{v-v^{-1}}$;

(5) For $i\not=j$, setting $b=1-a_{ij}$,
\begin{displaymath}
\sum_{k=0}^{b}(-1)^{k}E_i^{(k)}E_jE_i^{(b-k)}=0;
\end{displaymath}

(6) For $i\not=j$, setting $b=1-a_{ij}$,
\begin{displaymath}
\sum_{k=0}^{b}(-1)^{k}F_i^{(k)}F_jF_i^{(b-k)}=0.
\end{displaymath}
Here, ${K}_{i}=K_{h_i}$ and $E_i^{(n)}=E_i^n/[n]_{v}!$, $F_i^{(n)}=F_i^n/[n]_{v}!$.

Let $\mathbf{U}^+$ (resp. $\mathbf{U}^-$) be the subalgebra of $\mathbf{U}$ generated by the elements $E_i$ (resp. $F_i$) for $i\in I$, and let $\mathbf{U}^{0}$ be the subalgebra of $\mathbf{U}$ generated by $K_{\mu}$ for $\mu\in P^{\vee}$. We know that the quantized enveloping algebra has the triangular decomposition
\begin{displaymath}
\mathbf{U}\cong {\mathbf{U}^-}\otimes{\mathbf{U}^{0}}\otimes{\mathbf{U}^{+}}.
\end{displaymath}

We denote by $\bar{()}$ the unique automorphism of $\mathbf{U}$ as $\mathbb{Q}$-algebra given by
\begin{eqnarray*}
\bar{E}_i=E_i,\bar{F}_i=E_i,\bar{K}_{\mu}=K_{-\mu},\,\,
\textrm{for}\,\,i\in I,\mu\in p^{\vee},\\ \overline{fx}=\bar{f}\bar{x},\,  \textrm{for}\,  f\in \mathbb{Q}(v),x\in \mathbf{U}
\end{eqnarray*}where $\bar{f}(v)=f(v^{-1}).$

Let $\mathbf{f}$ be the associative algebra defined by Lusztig in [13]. Then $\mathbf{f}$ is generated by $\theta_i(i\in I)$ subject to  the above relation (6) with replacing $F_i$ by $\theta_i.$ Hence there are well-defined $\mathbb{Q}(v)$-algebra homomorphisms $\mathbf{f}\rightarrow \mathbf{U}(x\mapsto x^+)$ and $\mathbf{f}\rightarrow\mathbf{ U}(x\mapsto x^-)$ with image $\mathbf{U}^+$ and $\mathbf{U}^-$  respectively, where $E_i=\theta_i^+$ and $F_i=\theta_i^-$. The $\bar{()}$-involution of $\mathbf{U}$ induces a $\bar{()}$-involution of $\mathbf{f}$.

Let $\mathcal{A}=\mathbb{Q}[v,v^{-1}]$ and $\mathcal{Z}=\mathbb{Z}[v,v^{-1}]$. Denote by $\mathbf{U}^{\pm}_{\mathcal{Z}}$ the $\mathcal{Z}$-subalgebra of $\mathbf{U}^{\pm}$ generated by $E_{i}^{(s)}$ and $F_{i}^{(s)}$ for $i\in I$ and $s\in\mathbb{Z}$ respectively.
Also, we denote by $\mathbf{U}_{\mathcal{Z}}$ the $\mathcal{Z}$-subalgebra of $\mathbf{U}$ generated by $E_{i}^{(s)}$, $F_{i}^{(s)}$ and $K_{\mu}$ for $i\in I$, $s\in\mathbb{Z}$ and $\mu\in P^{\vee}$. Let $\mathbf{U}^{\pm}_{\mathcal{A}}=\mathbf{U}^{\pm}_{\mathcal{Z}}\otimes_{\mathcal{Z}}\mathcal{A}$ and $\mathbf{U}_{\mathcal{A}}=\mathbf{U}_{\mathcal{Z}}\otimes_{\mathcal{Z}}\mathcal{A}$. Similarly, let $\mathbf{f}_{\mathcal{Z}}$ be the $\mathcal{Z}$-subalgebra of $\mathbf{f}$ generated by $\theta_{i}^{(s)}$ for $i\in I$ and $s\in\mathbb{Z}$ and  $\mathbf{f}_{\mathcal{A}}=\mathbf{f}_{\mathcal{Z}}\otimes_{\mathcal{Z}}\mathcal{A}$.

In [10], [11] and [13], Lusztig defined the canonical basis of $\mathbf{f}$. We denote it by $\mathbf{B}$.

\subsec{2.2\quad Modified quantized enveloping algebra}

Let us recall the definition of the modified form $\dot{\mathbf{U}}$ of $\mathbf{U}$ in [13].

If $\lambda',\lambda''\in P$, we set
\begin{displaymath}
_{\lambda'}\mathbf{U}_{\lambda''}=\mathbf{U}/\left(\sum_{\mu\in P^{\vee}}(K_{\mu}-q^{\lambda'(\mu)})\mathbf{U}+(\sum_{\mu\in P^{\vee}}\mathbf{U}(K_{\mu}-q^{\lambda''(\mu)})\right).
\end{displaymath}
Let $\pi_{\lambda',\lambda''}:\mathbf{U}\rightarrow _{\lambda'}\mathbf{U}_{\lambda''}$ be the canonical projective and
\begin{displaymath}
\dot{\mathbf{U}}=\bigoplus_{\lambda',\lambda''\in P} {_{\lambda'}\mathbf{U}_{\lambda''}}.
\end{displaymath}

Consider the weight space decomposition $\mathbf{U}=\oplus_{\beta}\mathbf{U}(\beta)$, where $\beta$ through $Q$ and
$\mathbf{U}(\beta)=\{x\in\mathbf{U}|K_{\mu}xK_{\mu}^{-1}=v^{\beta(\mu)}x,\textrm{for $\mu\in P^{\vee}$}\}$.
The image of summands $\mathbf{U}(\beta)$ under $\pi_{\lambda',\lambda''}$ form the weight space decomposition $_{\lambda'}\mathbf{U}_{\lambda''}=\oplus_{\beta}{_{\lambda'}\mathbf{U}_{\lambda''}}(\beta)$.
Note that $_{\lambda'}\mathbf{U}_{\lambda''}(\beta)=0$ unless $\lambda'-\lambda''=\beta$.

There is a natural associative $\mathbb{Q}(v)$-algebra structure on $\dot{\mathbf{U}}$ inherited from that of $\mathbf{U}$. It is defined as follows: for any $\lambda'_1,\lambda''_1,\lambda'_2,\lambda''_2$, $\beta_1,\beta_2\in Q$ such that $\lambda'_1-\lambda''_1=\beta_1,\lambda'_2-\lambda''_2=\beta_2$ and any $x\in \mathbf{U}(\beta_1),y\in \mathbf{U}(\beta_2)$,
\begin{displaymath}
\pi_{\lambda'_1,\lambda''_1}(x)\pi_{\lambda'_2,\lambda''_2}(y)=
\left\{
  \begin{array}{cc}
    \pi_{\lambda'_1,\lambda''_2}(xy) & \textrm{if}\,\,\,\lambda''_1=\lambda'_2\\
     0 & \textrm{otherwise}
  \end{array}
\right..
\end{displaymath}

Let $\mathbf{1}_{\lambda}=\pi_{\lambda,\lambda}(\mathbf{1})$, where $\mathbf{1}$ is the unit element of $\mathbf{U}$. Then they satisfy $\mathbf{1}_{\lambda}\mathbf{1}_{\lambda'}=\delta_{\lambda,\lambda'}\mathbf{1}_{\lambda}$. In general, there is no unit element in the algebra $\dot{\mathbf{U}}$. However the family $(\mathbf{1}_{\lambda})_{\lambda\in P}$ can be regarded locally as the unit element in $\dot{\mathbf{U}}$.

Note that ${_{\lambda'}\mathbf{U}_{\lambda''}}=\mathbf{1}_{\lambda'}\dot{\mathbf{U}}\mathbf{1}_{\lambda''}$. We define $\dot{\mathbf{U}}\mathbf{1}_{\lambda}=\oplus_{\lambda'\in{P}}\mathbf{1}_{\lambda'}\dot{\mathbf{U}}\mathbf{1}_{\lambda}$. Then $\dot{\mathbf{U}}=\oplus_{\lambda\in{P}}\dot{\mathbf{U}}\mathbf{1}_{\lambda}$.

The $\mathbb{Q}$-algebra automorphism $\bar{()}:\mathbf{U}\rightarrow \mathbf{U}$ induces, for each $\lambda',\lambda''$, a linear isomorphism ${_{\lambda'}\mathbf{U}_{\lambda''}}\rightarrow{_{\lambda'}\mathbf{U}_{\lambda''}}$. Taking direct sums, we obtain an algebra automorphism $\bar{()}:\dot{\mathbf{U}}\rightarrow\dot{\mathbf{U}}$ which maps each $\mathbf{1}_{\lambda}$ to itself.

By Lusztig ([13], 23.2.1), we know that the elements $b^+b'^-\mathbf{1}_{\lambda}$ for $b,b'\in \mathbf{B}$ form a basis of the $\mathbb{Q}(v)$-vector space $\dot{\mathbf{U}}\mathbf{1}_{\lambda}$. Hence, the elements $b^+b'^-\mathbf{1}_{\lambda}$ for $b,b'\in \mathbf{B}$, $\lambda\in P$ form a basis of the $\mathbb{Q}(v)$-vector space $\dot{\mathbf{U}}$. This induces the triangular decomposition of $\dot{\mathbf{U}}$.

We denote by $\dot{\mathbf{U}}_{\mathcal{Z}}$ the subalgebra generated by the elements $E_{i}^{(n)}\mathbf{1}_{\lambda}$ and $F_{i}^{(n)}\mathbf{1}_{\lambda}$ over $\mathcal{Z}$ for all $i\in I$, $n\geq0$ and $\lambda\in P$. Then the elements $b^+b'^-\mathbf{1}_{\lambda}$ for $b,b'\in \mathbf{B}$, $\lambda\in P$ form an $\mathcal{Z}$-basis of $\dot{\mathbf{U}}_{\mathcal{Z}}$.

Lusztig ([13]) also defines the canonical basis of $\dot{\mathbf{U}}$ as a limit of the canonical basis of $L(\lambda)\otimes L^*(\mu)$ where $L(\lambda)$ is a highest weight module of $\mathbf{U}$ and $L^*(\mu)$ is a lowest weight module of $\mathbf{U}$. As the notation in [13], $\dot{\mathbf{B}}=\{b\diamondsuit_{\zeta} b'|b,b'\in{\mathbf{B}},\zeta\in P\}$ is the canonical basis of $\dot{\mathbf{U}}$. Note that $\{b\diamondsuit_{\lambda} b'|b,b'\in{\mathbf{B}}\}$ is a $\mathcal{Z}$-basis of $\dot{\mathbf{U}}_{\mathcal{Z}}\mathbf{1}_{\lambda}$. According to the proof of Theorem 25.2.1 in [13], we know that
\begin{displaymath}
b\diamondsuit_{\lambda} b'\equiv b^+b'^-\mathbf{1}_{\lambda}\,\,\,\textrm{mod}\,P(tr|b|-1,tr|b'|-1).
\end{displaymath}
Here $P(tr|b|-1,tr|b'|-1)$ is the $\mathbb{Q}(v)$-submodule of $\dot{\mathbf{U}}$ spanned by the set
\begin{displaymath}
\{b_1^+b_2^-1_{\zeta}|b_1,b_2\in \mathbf{B}\,\,\,\textrm{such that $tr|b_1|\leq tr|b|-1$, $tr|b_2|\leq tr|b'|-1$ and $|b_1|-|b_2|=|b|-|b'|$}\},
\end{displaymath}
where $|b|$ is the weight of $b$ and $tr\mu=\sum a_i$ for $\mu=\sum a_i\alpha_i$.

\subsec{2.3\quad Ringel-Hall algebra}

In this subsection, we recall the definition of Ringel-Hall algebras, following the notations in [9].

A quiver $Q=(I,H,s,t)$ consists of a vertex set $I$, an arrow set $H$, and two maps $s,t:H\rightarrow I$ such that an arrow $\rho\in H$ starts at $s(\rho)$ and terminates at $t(\rho)$.

We denote by $\mathbb{F}_q$ a finite field with $q$ elements and $\Lambda=\mathbb{F}_q(Q)$ the path algebra of $Q$ over $\mathbb{F}_q$. By mod-$\Lambda$ we denote the category of all finite dimension left $\Lambda$-modules. It is well-known that mod-$\Lambda$ is equivalent to the category of finite dimension representations of $Q$ over $\mathbb{F}_q$. We shall identify $\Lambda$-modules with representations of $Q$.

Given three modules $L,M$ and $N$ in mod-$\Lambda$, let $g^{L}_{MN}$ denote the number of $\Lambda$-submodules $W$ of $L$ such that $W\simeq N$ and $L/W\simeq N$ in mod-$\Lambda$. Let $v_q=\sqrt{q}\in \mathbb{C}$, $\mathcal{P}$ be the set of isomorphism classes of finite dimension nilpotent $\Lambda$-modules and ind$(\mathcal{P})$ be the set of isomorphism classes of indecomposable finite dimension nilpotent $\Lambda$-modules. The Ringel-Hall algebra $\mathcal{H}_q(\Lambda)$ of $\Lambda$ is by definition the $\mathbb{Q}(v_q)$-space with basis $\{u_{[M]}|[M]\in\mathcal{P}\}$ whose multiplication is given by
\begin{displaymath}
u_{[M]}u_{[N]}=\sum_{[L]\in\mathcal{P}}g^{L}_{MN}u_{[L]}.
\end{displaymath}
It is easily seen that $\mathcal{H}_q(\Lambda)$ is associative $\mathbb{Q}(v_q)$-algebra with unit $u_{0}$, where $0$ denotes the zero modules.

The set of isomorphism classes of (nilpotent) simple $\Lambda$-modules is naturally indexed by the set $I$ of vertices of $Q$. Then the Grothendieck group $G(\Lambda)$ of mod-$\Lambda$ is the free abelian group $\mathbb{Z}I$. For each nilpotent $\Lambda$-module $M$, the dimension vector $\underline{\dim}M=\sum_{i\in I}(\dim M_i)i$ is an element in $G(\Lambda)$. The Ringel-Hall algebra $\mathcal{H}_q(\Lambda)$ is a $\mathbb{N}I$-graded algebra by dimension vectors of modules.

The Euler form $\langle-,-\rangle$ on $G(\Lambda)=\mathbb{Z}I$ is defined by
\begin{displaymath}
\langle\alpha,\beta\rangle=\sum_{i\in I}a_ib_i-\sum_{\rho\in H}a_{s(\rho)}b_{t(\rho)}
\end{displaymath}
for $\alpha=\sum_{i\in I}a_{i}i$ and $\beta=\sum_{i\in I}b_{i}i$ in $\mathbb{Z}I$. For any nilpotent $\Lambda$-modules $M$ and $N$ one has
\begin{displaymath}
\langle\underline{\dim}M,\underline{\dim}N\rangle=\dim_{\mathbb{F}_q}\textrm{Hom}_{\Lambda}(M,N)-\dim_{\mathbb{F}_q}\textrm{Ext}_{\Lambda}(M,N).
\end{displaymath}
The symmetric Euler form is defined as $(\alpha,\beta)=\langle\alpha,\beta\rangle+\langle\beta,\alpha\rangle$ for $\alpha,\beta\in \mathbb{Z}I$. This gives rise to a symmetric generalized Cartan matrix $A=(a_{ij})_{i,j\in I}$ with $a_{ij}=(i,j)$. $A$ is independent of the field $\mathbb{F}_q$ and the orientation of $Q$.

The twisted Ringel-Hall algebra $\mathcal{H}^{\ast}_q(\Lambda)$ is defined as follows. Set $\mathcal{H}^{\ast}_q(\Lambda)=\mathcal{H}_q(\Lambda)$ as $\mathbb{Q}(v_q)$-vector space and define the multiplication by
\begin{displaymath}
u_{[M]}\ast u_{[N]}=v_q^{\langle\underline{\dim}M,\underline{\dim}N\rangle}\sum_{[L]\in\mathcal{P}}g^{L}_{MN}u_{[L]}.
\end{displaymath}
The composition algebra $\mathcal{C}^{\ast}_q(\Lambda)$ is a subalgebra of $\mathcal{H}^{\ast}_q(\Lambda)$ generated by $u_{[S_i]}$, $i\in I$, where $S_i$ is the nilpotent simple module corresponding to $i\in I$.
For any $\Lambda$-module $M$, we denote $\langle M\rangle=v_q^{-\dim M+\dim \textrm{End}_{\Lambda}(M)}u_{[M]}$. Note that $\{\langle M\rangle|M\in\mathcal{P}\}$ is a $\mathbb{Q}(v_q)$-basis of $\mathcal{H}^{\ast}_q(\Lambda)$.

Let $\mathcal{A}_{v_q}$ be the subring of $\mathbb{C}$ as the image of $\mathcal{A}$ under the map $\mathcal{A}\rightarrow \mathbb{C}$ sending $v$ to $v_q$. Let $\mathcal{C}_q^{\ast}(\Lambda)_{\mathcal{A}_{v_q}}$ be the $\mathcal{A}_{v_q}$-subalgebra of $\mathcal{H}_q^{\ast}(\Lambda)$ generated by $u_{[S_i]}^{(\ast m)}=\frac{u_{[S_i]}^{\ast m}}{[m]_{v_q}!}$, where $i\in I$ and $S_i$ is simple $\Lambda$-module.


Then we consider the generic form of Ringel-Hall algebra. Let $Q$ be a finite quiver and $\mathcal{H}_q^{\ast}(kQ)$ the twisted Ringel-Hall algebra of the path algebra $kQ$. Let $\mathcal{K}$ be a set of finite fields $k$ such that the set $\{q_k=|k||k\in\mathcal{K}\}$ is infinite. Let ${R}$ be an integral domain containing $\mathbb{Q}$ and an element $v_{q_k}$ such that $v_{q_k}^2=q_k$ for each $k\in\mathcal{K}$. For each $k\in\mathcal{K}$, we consider the composition algebra $\mathcal{C}^{\ast}_q(kQ)$ which is the ${R}$-algebra generated by the elements $u_{[S_i]}(k)$ of the Ringel-Hall algebra $\mathcal{H}_q^{\ast}(kQ)$. Consider the direct product
\begin{displaymath}
\mathcal{H}^{\ast}(Q)=\prod_{k\in\mathcal{K}}\mathcal{H}_q^{\ast}(kQ)
\end{displaymath}
and the elements $v=(v_{q_k})_{k\in\mathcal{K}}$, $v^{-1}=(v_{q_k}^{-1})_{k\in\mathcal{K}}$ and $u_{[S_i]}=(u_{[S_i]}(k))_{k\in\mathcal{K}}$. By $\mathcal{C}^{\ast}(Q)_{\mathcal{A}}$ we denote the subalgebra of $\mathcal{H}^{\ast}(Q)$ generated by $v$, $v^{-1}$ and $u_{[S_i]}$ over  $\mathbb{Q}$. We may regard it as the $\mathcal{A}$-algebra generated by $u_{[S_i]}$ where $v$ is considered as an indeterminate. Finally, we define the $\mathbb{Q}(v)$-algebra $\mathcal{C}^{\ast}(Q)=\mathbb{Q}(v)\otimes\mathcal{C}^{\ast}(Q)_{\mathcal{A}}$, called the generic twisted composition algebra of type $Q$.


\th {Remark 2.1.}\ {\it If $Q$ is a Dynkin quiver, then the generic composition algebra of $Q$ can be defined directly using Hall polynomials.}

Then we have the following well-known result of Green and Ringel ([4][18]).

\th {Theorem 2.1.}\ {\it Let $Q$ be a connected quiver, $A$ be the associated  generalized Cartan matrix, and $\mathbf{f}$ be the Lusztig's algebra of type $A$. Then the correspondence $u_{[S_i]}\mapsto\theta_i$, $i\in I$ induces an algebra isomorphism from $\mathcal{C}^{\ast}(Q)$ to $\mathbf{f}$.
}

\subsec{2.4\quad Root categories}


Given a Hom-finite, Krull-Schmidt triangulated category $\mathcal{C}$, we consider the Grothendieck group $G(\mathcal{C})$. That is, $G(\mathcal{C})$ is the quotient of a free abelian group with a basis $\{[M]|M\in \mathcal{C}\}$, indexed by the isomorphism classes of all objects in $\mathcal{C}$, subject to the relations $[X]-[Y]+[Z]$ provided there exist triangles of form $X\rightarrow Y\rightarrow Z\rightarrow TX$. For any $M\in \mathcal{C}$, we denote by $\underline{\dim} M$ the canonical image of $M$ in $G(\mathcal{C})$, called the dimension vector of $M$.

A triangulated category is called a $2$-period triangulated category if the translation $T$ satisfies $T^2\simeq $id.

Let $k$ be a field. Given a finite dimensional hereditary $k$-algebra $\Lambda$, we have the derived category $D^b(\Lambda)$ of mod-$\Lambda$ , which is a triangulated category with the translation $T$, obtained from the category of bounded complexes over mod-$\Lambda$ by localizing with respect to the set of all quasi-isomorphisms ([5][6]). For our purpose, we should consider the orbit category $\mathcal{R}(\Lambda)=D^b(\Lambda)/T^2$ of $D^b(\Lambda)$ under the equivalent functor $T^2$. Let $F:D^b(\Lambda)\rightarrow\mathcal{R}(\Lambda)$ be the canonical functor. The translation $T$ of $D^b(\Lambda)$ induces an equivalent functor of $\mathcal{R}(\Lambda)$ of order $2$, which is still denoted by $T$. By [14], we know that $\mathcal{R}(\Lambda)$ is also a triangulated category with $T$ as the translation and the covering functor $F:D^b(\Lambda)\rightarrow\mathcal{R}(\Lambda)$ sends each triangle in $D^b(\Lambda)$ to a triangle in $\mathcal{R}(\Lambda)$.

It is clear that the root category $\mathcal{R}=\mathcal{R}(\Lambda)$ is a $2$-period triangulated category.

Let $S_1,S_2,\cdots,S_n$ be the set of all non-isomorphism simple $\Lambda$-modules. Then we can consider a matrix $A=(a_{ij})_{n\times n}$ given by
\begin{displaymath}
\begin{array}{c}
  a_{ij}=\dim(\textrm{Hom}_{\mathcal{R}}(S_i,S_j))-\dim(\textrm{Hom}_{\mathcal{R}}(S_i,TS_j)) \\
  +\dim(\textrm{Hom}_{\mathcal{R}}(S_j,S_i))+\dim(\textrm{Hom}_{\mathcal{R}}(S_j,TS_i)).
\end{array}
\end{displaymath}
This matrix is a symmetric generalized Cartan matrix. The Cartan matrix $A$ only depend on the underlying graph of $Q$ and is independent of the the choice of the field and of the orientation of $Q$. Actually, the Dynkin graph of $A$ is just a complete slice of the quiver $Q$. So we can call $A$ the Cartan matrix of the root category $\mathcal{R}$ ([15]).

Let $Q$ be a connected quiver, $\mathcal{R}(Q)=D^b(kQ)/T^2$ be the root category. Let $\tilde{\mathcal{P}}$ be the set of isomorphism classes of the objects in $\mathcal{R}(Q)$ and ind$(\tilde{\mathcal{P}})$ be the set of isomorphism classes of the indecomposable objects in $\mathcal{R}(Q)$.
Note that mod-$kQ$ can be embedding into $\mathcal{R}(Q)$ as a full subcategory. Then ind$(\tilde{\mathcal{P}})=\textrm{ind($\mathcal{P}$)}\dot{\cup}\textrm{ind$(T(\mathcal{P}))$}$ where $\dot{\cup}$ is disjoint union.

\sec{3\quad Finite type}

\subsec{3.1\quad PBW-type basis of $\mathbf{U}^+$}

In this section, we consider a connected quiver $Q$ of Dynkin type.  We first consider the category of representations of $Q$ over some finite field $k$. Let $\Lambda=kQ$. We denote by $\Phi^+(\Phi^-)$ the set of positive roots (negative roots) of the Dynkin quiver Q. By the Gabriel's Theorem, we know that $\underline{\dim}$ induces an bijection between the set ind$(\mathcal{P})$ and the set $\Phi^+$. Given a positive root $\alpha$, we denote by $M(\alpha)$ the corresponding indecomposable representation of $Q$.

Since $Q$ is representation-directed, we can define a total order on the set
\begin{displaymath}
\Phi^+=\{\alpha_1,\alpha_2,\cdots,\alpha_n\}
\end{displaymath}
with
\begin{displaymath}
\{M(\alpha_1),M(\alpha_2),\cdots,M(\alpha_n)\}
\end{displaymath}
being the corresponding indecomposable $\Lambda$-modules such that
\begin{displaymath}
\textrm{Hom}(M(\alpha_i),M(\alpha_j))\not=0\Rightarrow i\leq j.
\end{displaymath}

We denote by $\mathbb{N}^{\textrm{ind$(\mathcal{P})$}}$ the set of all functions $\mathbf{a}:\Phi^+\rightarrow\mathbb{N}$. Each $\mathbf{a}\in\mathbb{N}^{\textrm{ind$(\mathcal{P})$}}$ defines a representation
\begin{displaymath}
M(\mathbf{a})=\bigoplus_{\alpha\in\Phi^+}\mathbf{a}(\alpha)M(\alpha)
\end{displaymath}
and any representation is isomorphic to one of the form.

Since the Hall polynomials exist in this case, we can consider the generic form $\mathcal{C}^{\ast}_{\mathcal{A}}(Q)$ of the Ringel-Hall algebra.

By [17], we have

\th {Proposition 3.1.}\ {\it The set $\{\langle M(\mathbf{a})\rangle|\mathbf{a}\in\mathbb{N}^{\textrm{ind$(\mathcal{P})$}}\}$ is a $\mathcal{A}$-basis of $\mathcal{C}^{\ast}_{\mathcal{A}}(Q)$.
}



\subsec{3.2\quad PBW-type basis of $\dot{\mathbf{U}}\mathbf{1}_{\lambda}$}

Consider a root category $\mathcal{R}(Q)$ over some finite field $k$. Remember that ind$(\tilde{\mathcal{P}})$ is the set of isomorphism classes of indecomposable objects in  $\mathcal{R}(Q)$. Let $\Phi=\{\underline{\dim}(M)|M\in\textrm{ind$(\tilde{\mathcal{P}})$}\}$. Then $\Phi$ is the root system of the corresponding Lie algebra and there is a bijection between the set ind$(\tilde{\mathcal{P}})$ and the set $\Phi$ by the Gabriel's Theorem. We have $\Phi=\Phi^+\dot{\cup}\Phi^-$. Given an element $\alpha\in\Phi$, we also denote by $M(\alpha)$ the corresponding object in $\mathcal{R}(Q)$. Let $\mathbb{N}^{\textrm{ind$(\tilde{\mathcal{P}})$}}$ be the set of all functions $\tilde{\mathbf{a}}:\Phi\rightarrow\mathbb{N}$. Each $\tilde{\mathbf{a}}\in\mathbb{N}^{\textrm{ind$(\tilde{\mathcal{P}})$}}$ defines an object
\begin{displaymath}
M(\tilde{\mathbf{a}})=\bigoplus_{\alpha\in\Phi}\tilde{\mathbf{a}}(\alpha)M(\alpha)
\end{displaymath}
and any object in our $\mathcal{R}(Q)$ is isomorphic to one of the form.

Note that the category $\mathcal{R}(Q)$, so the set $\mathbb{N}^{\textrm{ind$(\tilde{\mathcal{P}})$}}$, depend only on the underlying graph of Q.
If $Q'$ is another quiver such that $D^b(kQ)\simeq D^b(kQ')$, they give the same $\mathbb{N}^{\textrm{ind$(\tilde{\mathcal{P}})$}}$.

Given any symmetric generalized Cartan matrix $A=(a_{ij})_{n\times n}$ of finite type, we consider a quiver $Q$, the quantum enveloping algebra $\mathbf{U}$ and the modified enveloping algebra $\dot{\mathbf{U}}$ corresponding to the Cartan matrix $A=(a_{ij})_{n\times n}$.

Remember that mod-$kQ$ can be embedding into $\mathcal{R}(Q)$ as a full subcategory. Then ind$(\tilde{\mathcal{P}})=\textrm{ind($\mathcal{P}$)}\dot{\cup}\textrm{ind$(T(\mathcal{P}))$}$. For $\tilde{\mathbf{a}}\in\mathbb{N}^{\textrm{ind$(\tilde{\mathcal{P}})$}}$, let $\mathbf{a}_1=\tilde{\mathbf{a}}|_{\textrm{ind$(\mathcal{P})$}}$ and $\mathbf{a}_2=\tilde{\mathbf{a}}|_{\textrm{ind$(T(\mathcal{P}))$}}$, and we can denote by $\tilde{\mathbf{a}}=(\mathbf{a}_1,\mathbf{a}_2)$. We identify $\mathcal{C}^{\ast}(Q)$ with $\mathbf{f}$ by the correspondence between $u_{[S_i]}$ and $\theta_i$. So the set
\begin{displaymath}
\{\langle{M}(\tilde{\mathbf{a}})\rangle_{\lambda}=\langle{M}({\mathbf{a}_1})\rangle^+\cdot\langle{M}({\mathbf{a}_2})\rangle^-\mathbf{1}_{\lambda}
|\tilde{\mathbf{a}}\in\mathbb{N}^{\textrm{ind$(\tilde{\mathcal{P}})$}}\}
\end{displaymath}
can be regarded as elements in $\dot{\mathbf{U}}\mathbf{1}_{\lambda}$.


We have the following proposition

\th {Proposition 3.2.}\ {\it The set $\{\langle{M}(\tilde{\mathbf{a}})\rangle_{\lambda}|\tilde{\mathbf{a}}\in \mathbb{N}^{\textrm{ind$(\tilde{\mathcal{P}})$}}\}$ is a PBW-type basis of $\dot{\mathbf{U}}\mathbf{1}_{\lambda}$.}

\pf{Proof.} In ([13], 23.2.1), Lusztig points out that $\dot{\mathbf{U}}$ is a free $\mathbf{f}\otimes\mathbf{f}^{\textrm{opp}}$-module with basis $(\mathbf{1}_{\lambda})_{\lambda\in P}$. So the set
\begin{displaymath}
\{\langle{M}({\mathbf{a}_1})\rangle^+\cdot\langle{M}({\mathbf{a}_2})\rangle^-\mathbf{1}_{\lambda}
|\mathbf{a}_1\in{\textrm{ind$(\mathcal{P})$}},\mathbf{a}_2\in{\textrm{ind$(T(\mathcal{P}))$}}\}
\end{displaymath}
is a PBW-type basis of $\dot{\mathbf{U}}\mathbf{1}_{\lambda}$. By $\tilde{\mathbf{a}}=(\mathbf{a}_1,\mathbf{a}_2)$, we have the proposition.

\qed

We denote  by $B_{Q}(\dot{\mathbf{U}}\mathbf{1}_{\lambda})$ the PBW-type basis $\{\langle{M}(\tilde{\mathbf{a}})\rangle_{\lambda}|\tilde{\mathbf{a}}\in \mathbb{N}^{\textrm{ind$(\tilde{\mathcal{P}})$}}\}$. Note that the PBW-type basis depend on the embedding of mod-$kQ$ into $\mathcal{R}(Q)$.

\subsec{3.3\quad A bar invariant basis of $\dot{\mathbf{U}}\mathbf{1}_{\lambda}$}

Let $Q$ and $\mathcal{R}(Q)$ as before. Remember that $\Phi^+=\{\alpha_1,\dots,\alpha_n\}$. For $\mathbf{a},\mathbf{b}:\Phi^+\rightarrow\mathbb{N}$, we define $\mathbf{b}\prec\mathbf{a}$ if and only if there exists some $1\leq j\leq n$ such that $\mathbf{b}(\alpha_i)=\mathbf{a}(\alpha_i)$ for all $i<j$ and $\mathbf{b}(\alpha_j)>\mathbf{a}(\alpha_j)$.
For $\tilde{\mathbf{a}},\tilde{\mathbf{b}}:\Phi\rightarrow\mathbb{N}$, we define $\tilde{\mathbf{a}}\prec\tilde{\mathbf{b}}$ if and only if $\mathbf{a}_1\preceq\mathbf{b}_1$ and $\mathbf{a}_2\preceq\mathbf{b}_2$ but $\tilde{\mathbf{a}}\neq\tilde{\mathbf{b}}$, where $\tilde{\mathbf{a}}=(\mathbf{a}_1,\mathbf{a}_2)$ and $\tilde{\mathbf{b}}=(\mathbf{b}_1,\mathbf{b}_2)$.

Recall that for $\mathbf{c}:\Phi^+\rightarrow\mathbb{N}$, there exist a monomials $w_{\ast}(\mathbf{c})$ on Chevalley generators $u_{[S_i]}$  satisfying
\begin{displaymath}
w_{\ast}(\mathbf{c})=\langle M(\mathbf{c})\rangle+\sum_{\mathbf{c}'\prec\mathbf{c}}a_{\mathbf{c}\mathbf{c}'}\langle M(\mathbf{c}')\rangle,
\end{displaymath}
where $a_{\mathbf{c}\mathbf{c}'}\in \mathcal{A}$ ([17]).

Let $a=(a_{\mathbf{c}\mathbf{c}'})$ be the transition matrix from $\{\langle M(\mathbf{c})\rangle|\mathbf{c}:\Phi^+\rightarrow\mathbb{N}\}$
to $\{w_{\ast}(\mathbf{c})|\mathbf{c}:\Phi^+\rightarrow\mathbb{N}\}$, where  $a_{\mathbf{c}\mathbf{c}}=1$ and $a_{\mathbf{c}\mathbf{c}'}=0$ unless $\mathbf{c}'\prec\mathbf{c}$. Note that $a$ is unipotent lower triangular matrix.

Let $\bar{a}$ be obtained from $a$ by applying the $\bar{()}$-involution to each elements of $a$.  Since $\overline{w_{\ast}(\mathbf{c})}=w_{\ast}(\mathbf{c})$, we have
\begin{displaymath}
w_{\ast}(\mathbf{c})=\overline{w_{\ast}(\mathbf{c})}=\sum_{\mathbf{c}'}\bar{a}_{\mathbf{c}\mathbf{c}'}\overline{\langle M(\mathbf{c}')\rangle},
\end{displaymath}
thus
\begin{displaymath}
\overline{\langle M(\mathbf{c})\rangle}=\sum_{\mathbf{c}'}{\bar{a}^{-1}}_{\mathbf{c}\mathbf{c}'}w_{\ast}(\mathbf{c}')
=\sum_{\mathbf{c}'}\sum_{\mathbf{c}''}{\bar{a}^{-1}}_{\mathbf{c}\mathbf{c}'}{a}_{\mathbf{c}'\mathbf{c}''}\langle M(\mathbf{c}'')\rangle.
\end{displaymath}
Let $h=\bar{a}^{-1}a$, then $h$ is again a unipotent lower triangular matrix, and $\bar{h}=h^{-1}$. There exists a unique unipotent lower triangular matrix $d=(d_{\mathbf{c}\mathbf{c}'})$ with off-diagonal entries in $v^{-1}\mathbb{Q}[v^{-1}]$, such that $d=\bar{d}h$. Then the canonical basis of $\mathbf{f}$ is
\begin{displaymath}
\mathcal{E}^{\mathbf{c}}=\langle M(\mathbf{c})\rangle+\sum_{\mathbf{c}'\prec\mathbf{c}}d_{\mathbf{c}\mathbf{c}'}\langle M(\mathbf{c}')\rangle,
\end{displaymath}
with $d_{\mathbf{c}\mathbf{c}'}\in v^{-1}\mathbb{Q}[v^{-1}]$ ([17]).

Similarly, we can get a bar-invariant basis of $\dot{\mathbf{U}}\mathbf{1}_{\lambda}$ from
\begin{displaymath}
B_{Q}(\dot{\mathbf{U}}\mathbf{1}_{\lambda})=\{\langle{M}({\mathbf{c}_1})\rangle^+\cdot\langle{M}({\mathbf{c}_2})\rangle^-\mathbf{1}_{\lambda}
|\tilde{\mathbf{c}}:\Phi\rightarrow\mathbb{N},\tilde{\mathbf{c}}=(\mathbf{c}_1,\mathbf{c}_2)\}
\end{displaymath}
and
\begin{displaymath}
\{w_{\ast}(\mathbf{c}_1)^+\cdot w_{\ast}(\mathbf{c}_2)^-\mathbf{1}_{\lambda}|\tilde{\mathbf{c}}:\Phi\rightarrow\mathbb{N},\tilde{\mathbf{c}}=(\mathbf{c}_1,\mathbf{c}_2)\}
\end{displaymath}
under the order $\prec$ on $\mathbb{N}^{\textrm{ind$(\tilde{\mathcal{P}})$}}$ defined above. We define $w_{\ast}(\tilde{\mathbf{c}})_{\lambda}=w_{\ast}(\mathbf{c}_1)^+\cdot{w_{\ast}(\mathbf{c}_2)}^-\mathbf{1}_{\lambda}$ where $\tilde{\mathbf{c}}=(\mathbf{c}_1,\mathbf{c}_2)$.

By the relation
\begin{displaymath}
w_{\ast}(\mathbf{c})=\langle M(\mathbf{c})\rangle+\sum_{\mathbf{c}'\prec\mathbf{c}}a_{\mathbf{c}\mathbf{c}'}\langle M(\mathbf{c}')\rangle,
\end{displaymath}
we have
\begin{displaymath}
w_{\ast}(\mathbf{c}_1)^+=\langle M(\mathbf{c}_1)\rangle^++\sum_{\mathbf{c}'_1\prec\mathbf{c}_1}a_{\mathbf{c}_1\mathbf{c}'_1}\langle M(\mathbf{c}'_1)\rangle^+
\end{displaymath}
and
\begin{displaymath}
w_{\ast}(\mathbf{c}_2)^-=\langle M(\mathbf{c}_2)\rangle^-+\sum_{\mathbf{c}'_2\prec\mathbf{c}_2}a_{\mathbf{c}_2\mathbf{c}'_2}\langle M(\mathbf{c}'_2)\rangle^-
\end{displaymath}
in $\mathbf{U}^{\pm}$ respectively.
Hence, we have
\begin{eqnarray*}
w_{\ast}(\tilde{\mathbf{c}})_{\lambda}&=&w_{\ast}(\mathbf{c}_1)^+\cdot w_{\ast}(\mathbf{c}_2)^-\mathbf{1}_{\lambda} \\
&=&(\langle M(\mathbf{c}_1)\rangle+\sum_{\mathbf{c}'_1\prec\mathbf{c}_1}a_{\mathbf{c}_1\mathbf{c}_1'}\langle M(\mathbf{c}'_1)\rangle)^+\cdot
(\langle M(\mathbf{c}_2)\rangle+\sum_{\mathbf{c}'_2\prec\mathbf{c}_2}a_{\mathbf{c}_2\mathbf{c}'_2}\langle M(\mathbf{c}'_2)\rangle)^-\mathbf{1}_{\lambda}\\
&=&\langle M(\mathbf{c}_1)\rangle^+\cdot\langle M(\mathbf{c}_2)\rangle^-\mathbf{1}_{\lambda}+
\langle M(\mathbf{c}_1)\rangle^+\cdot\sum_{\mathbf{c}'_2\prec\mathbf{c}_2}a_{\mathbf{c}_2\mathbf{c}'_2}\langle M(\mathbf{c}'_2)\rangle^-\mathbf{1}_{\lambda}+\\
&&\sum_{\mathbf{c}'_1\prec\mathbf{c}_1}a_{\mathbf{c}_1\mathbf{c}_1'}\langle M(\mathbf{c}'_1)\rangle^+\cdot\langle M(\mathbf{c}_2)\rangle^-\mathbf{1}_{\lambda}+
\sum_{\mathbf{c}'_1\prec\mathbf{c}_1}a_{\mathbf{c}_1\mathbf{c}_1'}\langle M(\mathbf{c}'_1)\rangle^+\cdot\sum_{\mathbf{c}'_2\prec\mathbf{c}_2}a_{\mathbf{c}_2\mathbf{c}'_2}\langle M(\mathbf{c}'_2)\rangle^-\mathbf{1}_{\lambda}\\
&=&\langle{M}(\tilde{\mathbf{c}})\rangle_{\lambda}
+\sum_{\tilde{\mathbf{c}}'\prec\tilde{\mathbf{c}}}\tilde{a}_{\tilde{\mathbf{c}}\tilde{\mathbf{c}}'}\langle{M}(\tilde{\mathbf{c}}')\rangle_{\lambda},
\end{eqnarray*}
where $\tilde{\mathbf{c}}=(\mathbf{c}_1,\mathbf{c}_2)$, $\tilde{\mathbf{c}}'=(\mathbf{c}'_1,\mathbf{c}'_2)$ and $\tilde{a}_{\tilde{\mathbf{c}}\tilde{\mathbf{c}}'}=a_{\mathbf{c}_1\mathbf{c}_1'}a_{\mathbf{c}_2\mathbf{c}'_2}$.

The same as above, let $\tilde{a}=(\tilde{a}_{\tilde{\mathbf{c}}\tilde{\mathbf{c}}'})$ be the transition matrix from $\{\langle M(\tilde{\mathbf{c}})\rangle_{\lambda}|\tilde{\mathbf{c}}:\Phi\rightarrow\mathbb{N}\}$
to $\{w_{\ast}(\tilde{\mathbf{c}})_{\lambda}|\tilde{\mathbf{c}}:\Phi\rightarrow\mathbb{N}\}$, where  $\tilde{a}_{\tilde{\mathbf{c}}\tilde{\mathbf{c}}}=1$ and $a_{\widetilde{\mathbf{c}}\tilde{\mathbf{c}}'}=0$ unless $\tilde{\mathbf{c}}'\prec\tilde{\mathbf{c}}$. Note that $\tilde{a}$ is unipotent lower triangular matrix with off-diagonal entries in $\mathcal{A}$.

Let $\bar{\tilde{a}}$ be obtained from $\tilde{a}$ by applying the $\bar{()}$-involution to each elements of $\tilde{a}$. Since $\overline{w_{\ast}(\tilde{\mathbf{c}})_{\lambda}}=w_{\ast}(\tilde{\mathbf{c}})_{\lambda}$, we have
\begin{displaymath}
w_{\ast}(\tilde{\mathbf{c}})_{\lambda}=\overline{w_{\ast}(\tilde{\mathbf{c}})_{\lambda}}=\sum_{\tilde{\mathbf{c}}'}\bar{a}_{\tilde{\mathbf{c}}\tilde{\mathbf{c}}'}\overline{\langle {M}(\tilde{\mathbf{c}}')\rangle_{\lambda}},
\end{displaymath}
thus
\begin{displaymath}
\overline{\langle{M}(\tilde{\mathbf{c}})\rangle_{\lambda}}=\sum_{\tilde{\mathbf{c}}'}\bar{\tilde{a}}^{-1}_{\tilde{\mathbf{c}}\tilde{\mathbf{c}}'}w_{\ast}(\tilde{\mathbf{c}}')_{\lambda}
=\sum_{\tilde{\mathbf{c}}'}\sum_{\tilde{\mathbf{c}}''}\bar{\tilde{a}}^{-1}_{\tilde{\mathbf{c}}\tilde{\mathbf{c}}'}{\tilde{a}}_{\tilde{\mathbf{c}}'\tilde{\mathbf{c}}''}\langle {M}(\tilde{\mathbf{c}}'')\rangle_{\lambda}.
\end{displaymath}

Let $\tilde{h}=\bar{\tilde{a}}^{-1}\tilde{a}$, then $\tilde{h}$ is again a unipotent lower triangular matrix, and $\bar{\tilde{h}}=\tilde{h}^{-1}$. There exists a unique unipotent lower triangular matrix $\tilde{d}=(\tilde{d}_{\tilde{\mathbf{c}}\tilde{\mathbf{c}}'})$ with off-diagonal entries in $v^{-1}\mathbb{Q}[v^{-1}]$, such that $\tilde{d}=\bar{\tilde{d}}\tilde{h}$. Then we can define a bar-invariant basis of $\dot{\mathbf{U}}\mathbf{1}_{\lambda}$
\begin{displaymath}
\mathcal{E}^{\tilde{\mathbf{c}}}_{\lambda}=\langle{M}(\tilde{\mathbf{c}})\rangle_{\lambda}+
\sum_{\tilde{\mathbf{c}}'\prec\tilde{\mathbf{c}}}\tilde{d}_{\tilde{\mathbf{c}}\tilde{\mathbf{c}}'}\langle{M}(\tilde{\mathbf{c}}')\rangle_{\lambda},
\end{displaymath}
with $\tilde{d}_{\tilde{\mathbf{c}}'\tilde{\mathbf{c}}}\in v^{-1}\mathbb{Q}[v^{-1}]$. We denoted by $\mathcal{B}_{Q}(\dot{\mathbf{U}}\mathbf{1}_{\lambda})$ the above basis.

\th {Theorem 3.1.}\ {\it $\mathcal{B}_{Q}(\dot{\mathbf{U}}\mathbf{1}_{\lambda})=\{\mathcal{E}^{\tilde{\mathbf{c}}}_{\lambda}|\tilde{\mathbf{c}}:\Phi\rightarrow\mathbb{N}\}
        =\{b^+b'^-\mathbf{1}_{\lambda}|b,b'\in \mathbf{B}\}$.}

We omit the proof of the above theorem. The proof of Theorem 3.1 is simple than Theorem 4.1 of affine case, which will be proved in next section.

\subsec{3.4\quad A parameterization of the canonical basis of $\dot{\mathbf{U}}\mathbf{1}_{\lambda}$}

Let $\dot{\mathbf{U}}=\oplus_{\lambda\in P}\dot{\mathbf{U}}\mathbf{1}_{\lambda}$ be the modified enveloping algebra corresponding to the quiver $Q$ and $\dot{\mathbf{B}}_{\lambda}$ is the canonical basis of $\dot{\mathbf{U}}\mathbf{1}_{\lambda}$.

\th {Theorem 3.2.}\ {\it We have a bijective map
\begin{displaymath}
\Psi_Q:\mathbb{N}^{\textrm{ind$(\tilde{\mathcal{P}})$}}\rightarrow{\dot{\mathbf{B}}_{\lambda}}
\end{displaymath}
given by
\begin{displaymath}
\tilde{\mathbf{c}}\mapsto\mathcal{E}^{\mathbf{c}_1}\diamondsuit_{\lambda}\mathcal{E}^{\mathbf{c}_2},
\end{displaymath}
which is the composition of the following two bijection
\begin{eqnarray*}
\mathbb{N}^{\textrm{ind$(\tilde{\mathcal{P}})$}}&\rightarrow&\mathcal{B}_{Q}(\dot{\mathbf{U}}\mathbf{1}_{\lambda})\\
\tilde{\mathbf{c}}&\mapsto&\mathcal{E}^{\tilde{\mathbf{c}}}_{\lambda},
\end{eqnarray*}
and
\begin{eqnarray*}
\mathcal{B}_{Q}(\dot{\mathbf{U}}\mathbf{1}_{\lambda})&\rightarrow&\dot{\mathbf{B}}_{\lambda}\\
b^+b'^-\mathbf{1}_{\lambda}&\mapsto&b\diamondsuit_{\lambda}b'.
\end{eqnarray*}
}

\pf{Proof.}\  The first bijection from $\mathbb{N}^{\textrm{ind$(\tilde{\mathcal{P}})$}}$ to $\mathcal{B}_{Q}(\dot{\mathbf{U}}\mathbf{1}_{\lambda})$ comes from our construction of $\mathcal{E}^{\tilde{\mathbf{c}}}_{\lambda}$ and the second bijection from $\mathcal{B}_{Q}(\dot{\mathbf{U}}\mathbf{1}_{\lambda})$ to $\dot{\mathbf{B}}_{\lambda}$ comes from Lusztig ([13], Theorem 25.2.1). By Theorem 3.1, $\mathcal{E}^{\tilde{\mathbf{c}}}_{\lambda}
=\mathcal{E}^{\mathbf{c}_1+}\mathcal{E}^{\mathbf{c}_2-}\mathbf{1}_{\lambda}$. So, we have the theorem.

\qed

For a non symmetric Lie algebra $\textbf{g}$, the similar result holds for $\dot{\mathbf{U}}(\textbf{g})$.

\sec{4\quad Affine type}

\subsec{4.1\quad PBW-type basis of $\mathbf{U}^+$}

We first recall the construction of the PBW-type basis in [9].

\subsec{4.1.1\quad The integral basis arising from the Kronecker quiver}\label{Kronecker}

Let $Q$ be the Kronecker quiver with $I=\{1,2\}$ and $H=\{\rho_1,\rho_2\}$ as follow
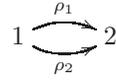
\begin{displaymath}
\xymatrix{
1\ar@/^/[r]^{\rho_1}\ar@/_/[r]_{\rho_2} & 2
}
\end{displaymath}

Let $\Lambda=\mathbb{F}_q(Q)$ be the path algebra.

The set of dimension vectors of indecomposable representations is
\begin{displaymath}
\Phi^{+}=\{(l+1,l),(m,m),(n,n+1)|l\geq 0,m\geq 1,n\geq 0\}.
\end{displaymath}
The dimension vectors $(l+1,l)$ and $(n,n+1)$ correspond to preprojective and preinjective indecomposable representations respectively.

Let $\mathcal{P}$ be the set of isomorphism classes of finite dimension $\Lambda$-modules, $\mathcal{H}_q$ (resp. $\mathcal{H}_q^{\ast}$) be the Ringel-Hall (resp. the twisted Ringel-Hall) algebra of $\Lambda$ over $\mathbb{Q}(v_q)$, where $v_q^2=q$.

Define
\begin{displaymath}
E_{(n+1,n)}=\langle u_{(n+1,n)}\rangle \,\,\textrm{and}\,\,E_{(n,n+1)}=\langle u_{(n,n+1)}\rangle.
\end{displaymath}
For $n\geq 1$, define
\begin{displaymath}
\tilde{E}_{n\delta}=E_{(n-1,n)}\ast E_{(1,0)}-v_q^{-2}E_{(1,0)}\ast E_{(n-1,n)}.
\end{displaymath}
Then, we can define by induction
\begin{displaymath}
E_{0\delta}=1,\,\,E_{k\delta}=\frac{1}{[k]}\sum_{s=1}^{k}v_{q}^{s-k}\tilde{E}_{s\delta}\ast E_{(k-s)\delta}\,\,\textrm{for $k\geq 1$}.
\end{displaymath}

Then we consider $\mathcal{C}^{\ast}_{\mathcal{Z}}$. Since $E_{k\delta}$, $E_{(m+1,m)}$ and $E_{(n,n+1)}$ can be defined in each $\mathcal{H}_q$, we can consider them as elements in $\prod_q\mathcal{H}_q^{\ast}$. We know that the set $\{E_{(m+1,m)},E_{k\delta},E_{(n,n+1)}|m,n\geq0,k\geq1\}$ is contained in $\mathcal{C}^{\ast}_{\mathcal{Z}}$.

Let $\mathbf{P}(m)$ be the set of all partition of $m$. For any partition
\begin{displaymath}
w=(w_1,w_2,\ldots,w_t)\in\mathbf{P}(m),
\end{displaymath}
we define
\begin{displaymath}
E_{w\delta}=E_{w_1\delta}\ast E_{w_2\delta}\ast\cdots\ast E_{w_t\delta}.
\end{displaymath}

\th {Proposition 4.1. [9]}\ {\it The set
\begin{displaymath}
\{\langle P\rangle\ast E_{w\delta}\ast \langle I\rangle\}
\end{displaymath}
where $P\in\mathcal{P}$ is preprojective, $w\in\mathbf{P}(m)$, $I\in\mathcal{P}$ is preinjective and $m\in\mathbb{Z}_{\geq 1}$, is a $\mathcal{Z}$-basis of $\mathcal{C}^{\ast}_{\mathcal{Z}}$.}

\subsec{4.1.2\quad The integral basis arising from a tube}\label{tube}

Let $\Delta=\Delta(n)$ be the cyclic quiver with vertex set $\Delta_{0}=\mathbb{Z}/n\mathbb{Z}$=\{1,2,\ldots,n\} and arrow set $\Delta_q=\{i\rightarrow i+1|i\in\Delta_{0}\}$ as follow
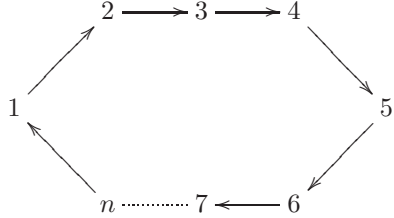
\begin{displaymath}
\xymatrix{
        &2\ar[r]&3\ar[r]&4\ar[dr]\\
1\ar[ur]&       &       &       &5\ar[dl]\\
        &n\ar[ul]&7\ar@{.}[l]&6\ar[l]
}
\end{displaymath}

We consider the category $\mathcal{T}=\mathcal{T}(n)$ of finite dimension nilpotent representations of $\Delta(n)$ over $\mathbb{F}_q$. Let $S_i,i\in \Delta_{0}$ be the irreducible objects in $\mathcal{T}(n)$ and $S_i[l]$ be the indecomposable objects in $\mathcal{T}(n)$ with top $S_i$ and length $l$. Note that $S_i[l]$ is independent of the choice of $q$. Let $\mathcal{P}$ be the set of isomorphism classes of objects in $\mathcal{T}(n)$. Denote by $\mathcal{H}$ (resp. $\mathcal{H}^{\ast}$) the Ringel-Hall algebra (resp. twisted Ringel-Hall algebra) of $\mathcal{T}(n)$. Because the Hall polynomials always exist in this case, we may regard them as generic form.

Let $\Pi$ be the set of $n$-tuples of partitions $\pi=(\pi^{(1)},\pi^{(2)},\ldots,\pi^{(n)})$ with each component $\pi^{(i)}=(\pi^{(i)}_1\geq\pi^{(i)}_2\geq\cdots)$ being a partition of integers. For each $\pi\in\Pi$, we define an object in $\mathcal{T}(n)$
\begin{displaymath}
M(\pi)=\bigoplus_{i\in\Delta_{0},j\geq 1}S_i[\pi_j^{(i)}].
\end{displaymath}
In this way we obtain a bijection between the set $\Pi$ and the set $\mathcal{P}$.

An $n$-tuple $\pi=(\pi^{(1)},\pi^{(2)},\ldots,\pi^{(n)})$ of partition in $\Pi$ is called aperiodic, if for each $l\geq1$ there is some $i=i(l)\in\Delta_{0}$ such that $\pi_j^{(i)}\neq l$ for all $j\geq 1$. By $\Pi^a$ we denote the set of aperiodic $n$-tuples of partitions. An object $M$ in $\mathcal{T}$ is called aperiodic if $M\simeq M(\pi)$ for some $\pi\in \Pi^a$. For any dimension vector $\alpha\in\mathbb{N}^n$, define
$\Pi_{\alpha}=\{\lambda\in\Pi|\underline{\dim}M(\lambda)=\alpha\}$ and $\Pi_{\alpha}^{a}=\Pi^{a}\cap\Pi_{\alpha}$.

Given any two modules $M,N$ in $\mathcal{T}$, there exists a unique (up to isomorphism) extension $L$ of $M$ by $N$ with minimal dimEnd$(L)$. The extension $L$ is called the generic extension of $M$ by $N$ and is denoted by $L=M\diamond N$.

Let $\Omega$ be the set of all words on the alphabet $\Delta_{0}$. For each $w=i_1 i_2\cdots i_m\in\Omega$, we set $M(w)=S_{i_1}\diamond S_{i_2}\diamond\cdots\diamond S_{i_m}$. Then there is a unique $\mathfrak{p}(w)=\pi\in\Pi$ such that $M(\pi)\simeq M(w)$. It has been proved in [16] that $\pi=\mathfrak{p}(w)\in\Pi^a$ and $\mathfrak{p}$ induces a surjective $\mathfrak{p}:\Omega\twoheadrightarrow \Pi^a$.

For each module $M$ in $\mathcal{T}$ and $s\geq 1$, we define by $sM$ the direct sum of $s$ copies of $M$.  For $w\in\Omega$, write $w$ in a tight form $w=j_1^{e_1}j_2^{e_2}\cdots j_t^{e_t}\in\Omega$ with $j_{r-1}\not=j_{r}$ for all $r$. We can get $\mu_r\in\Pi$ such that $M(\mu_r)=e_rS_{j_r}$. For any $\lambda\in \Pi_{\sum_{r=1}^{t}e_rj_r}$, write $g^{\lambda}_{w}$ for the Hall polynomial $g^{M(\lambda)}_{M(\mu_1),\ldots,M(\mu_t)}$. A word $w$ is caller distinguished if the Hall Polynomial $g^{\mathfrak{p}(w)}_{w}=1$. For any $\pi\in\Pi^a$, there exists a distinguished word $w_{\pi}=j_1^{e_1}j_2^{e_2}\cdots j_t^{e_t}\in\mathfrak{p}^{-1}(\pi)$ in tight form by [3]. From now on, we fix a distinguished word $w_{\pi} \in\mathfrak{p}^{-1}(\pi)$. Thus we have a section $\mathcal{D}=\{w_{\pi}|\pi\in\Pi^a\}$ of $\mathfrak{p}$ over $\Pi^a$. $\mathcal{D}$ is called a section of distinguished words in [3].

For each $w=j_1^{e_1}j_2^{e_2}\cdots j_t^{e_t}\in\Omega$ in tight form, define in $\mathcal{C}^{\ast}$ the monomial
\begin{displaymath}
\mathbf{m}^{(w)}=E_{j_i}^{\ast e_1}\ast\cdots\ast E_{j_t}^{\ast e_t}.
\end{displaymath}
Then define $E_{\pi}$ for all $\pi\in\Pi^a$ inductively by the following relation
\begin{displaymath}
E_{\pi}=\mathbf{m}^{(w_{\pi})}-\sum_{\lambda\prec\pi,\lambda\in\Pi^{a}_{\alpha}}v^{d}g^{\lambda}_{w_{\pi}}(v^2)E_{\lambda}
\end{displaymath}
and
\begin{displaymath}
E_{\pi}=\mathbf{m}^{(w_{\pi})}\,\,\textrm{if $\pi\in\Pi^{a}_{\alpha}$ is minimal},
\end{displaymath}
where $\alpha=\sum_{i=1}^{t} e_r j_r$, $d=-\dim M(\pi)+\dim \textrm{End} M(\pi)+\dim M(\lambda)-\dim \textrm{End} M(\lambda)$ and $\lambda\prec\mu\Leftrightarrow \dim \textrm{Hom}(M,M(\lambda))\leq \dim \textrm{Hom}(M,M(\mu))$ for all modules $M$ in $\mathcal{T}$. We know that $E_{\pi}\in \mathcal{C}^{\ast}$ for all $\lambda\in\Pi^a$.

\th {Proposition 4.2. [9]}\ {\it Let $\mathcal{D}=\{w_{\pi}|\pi\in\Pi^a\}$ be a section of distinguished words. Then both $\{\mathbf{m}^{(w_{\pi})}|\pi\in\Pi^a\}$ and $\{E_{\pi}|\pi\in\Pi^a\}$ are $\mathcal{Z}$-basis of $\mathcal{C}_{\mathcal{Z}}^{\ast}$.
And the transition matrix between these two basis is triangular with diagonal entries equal to $1$ and entries above the diagonal in $\mathbb{Z}[v,v^{-1}]$.
}

It has been proved in [3] that the basis $\{E_{\pi}|\pi\in\Pi^a\}$ is independent of the choice of the sections of distinguished words.

\subsec{4.1.3\quad The integral basis arising from preprojective and preinjective components}\label{prep and prei}

In this section, we consider a connected tame quiver $Q$ without oriented cycles. Let $\Lambda=\mathbb{F}_q(Q)$ be the path algebra. We denote by $Prep$ and $Prei$ the isomorphism classes of indecomposable preprojective and preinjective $\Lambda$-modules, which are independent of the choice of $q$. Let $\mathcal{H}_q$ (resp. $\mathcal{H}_q^{\ast}$) be the Ringel-Hall (resp. the twisted Ringel-Hall) algebra of $\Lambda$ over $\mathbb{Q}(v_q)$, where $v_q^2=q$.

Since $Prei$ is representation-directed, we can define a total order on the set
\begin{displaymath}
\Phi_{Prei}^{+}=\{\cdots,\beta_3,\beta_2,\beta_1\}
\end{displaymath}
of all positive real roots appearing in $Prei$ with
\begin{displaymath}
\{\cdots,M(\beta_3),M(\beta_2),M(\beta_1)\}
\end{displaymath}
being the corresponding indecomposable preinjective $\Lambda$-modules such that
\begin{displaymath}
\textrm{Hom}(M(\beta_i),M(\beta_j))\not=0\Rightarrow i\geq j.
\end{displaymath}

Similarly, since $Prep$ is representation-directed, we can define a total order on the set
\begin{displaymath}
\Phi_{Prep}^{+}=\{\alpha_1,\alpha_2,\alpha_3,\cdots\}
\end{displaymath}
of all positive real roots appearing in $Prep$ with
\begin{displaymath}
\{M(\alpha_1),M(\alpha_2),M(\alpha31),\cdots\}
\end{displaymath}
being the corresponding indecomposable preprojective $\Lambda$-modules such that
\begin{displaymath}
\textrm{Hom}(M(\alpha_i),M(\alpha_j))\not=0\Rightarrow i\leq j.
\end{displaymath}

We denote by $\mathbb{N}^{Prei}_{f}$ the set of all support-finite functions $\mathbf{b}:\Phi_{Prei}^{+}\rightarrow\mathbb{N}$. Each $\mathbf{b}\in\mathbb{N}^{Prei}_{f}$ defines a preinjective representation
\begin{displaymath}
M(\mathbf{b})=\bigoplus_{\beta_i\in\Phi_{Prei}^{+}}\mathbf{b}(\beta_i)M(\beta_i)
\end{displaymath}
and any preinjective representation is isomorphic to one of the form.

We denote by $\mathbb{N}^{Prep}_{f}$ the set of all support-finite functions $\mathbf{a}:\Phi_{Prep}^{+}\rightarrow\mathbb{N}$. Each $\mathbf{a}\in\mathbb{N}^{Prep}_{f}$ defines a preprojective representation
\begin{displaymath}
M(\mathbf{a})=\bigoplus_{\alpha_i\in\Phi_{Prep}^{+}}\mathbf{a}(\alpha_i)M(\alpha_i)
\end{displaymath}
and any preprojective representation is isomorphic to one of the form.

If $\mathbf{b},\mathbf{b}_1,\mathbf{b}_2\in \mathbb{N}^{Prei}_{f}$ (resp. $\mathbf{a},\mathbf{a}_1,\mathbf{a}_2\in \mathbb{N}^{Prep}_{f}$), the Hall polynomial $g^{M(\mathbf{b})}_{M(\mathbf{b}_1)M(\mathbf{b}_2)}$ (resp. $g^{M(\mathbf{a})}_{M(\mathbf{a}_1)M(\mathbf{a}_2)}$) always exists.

Then we consider $\mathcal{C}^{\ast}_{\mathcal{Z}}$, the generic form of twisted composition algebra of $Q$. We have $\langle M(\mathbf{b})\rangle\in\mathcal{C}^{\ast}_{\mathcal{Z}}$ (resp. $\langle M(\mathbf{a})\rangle\in\mathcal{C}^{\ast}_{\mathcal{Z}}$). We define $\mathcal{C}^{\ast}_{Prei}$ (resp. $\mathcal{C}^{\ast}_{Prep}$) to be the $\mathcal{Z}$-submodule of $\mathcal{C}^{\ast}_{\mathcal{Z}}$ generated by $\{\langle M(\mathbf{b})\rangle|\mathbf{b}\in\mathbb{N}^{Prei}_{f}\}$ (resp. $\{\langle M(\mathbf{a})\rangle|\mathbf{a}\in\mathbb{N}^{Prep}_{f}\}$).

\th {Proposition 4.3. [9]}\ {\it The $\mathcal{Z}$-submodule $\mathcal{C}^{\ast}_{Prei}$ (resp. $\mathcal{C}^{\ast}_{Prep}$) is an subalgebra of $\mathcal{C}^{\ast}_{\mathcal{Z}}$ and $\{\langle M(\mathbf{b})\rangle |\mathbf{b}\in\mathbb{N}^{Prei}_{f}\}$ (resp. $\{\langle M(\mathbf{a})\rangle|\mathbf{a}\in\mathbb{N}^{Prep}_{f}\}$) is a $\mathcal{Z}$-basis of $\mathcal{C}^{\ast}_{Prei}$ (resp. $\mathcal{C}^{\ast}_{Prep}$).}

\subsec{4.1.4\quad The integral basis for the generic composition algebras}

In this section, we still assume that $Q$ is connected tame quiver without oriented cycles. We first consider the embedding of the representation category of Kronecker quiver into the representation category of $Q$.

Let $e$ be a extending vertex of $Q$ and $\Lambda=\mathbb{F}_qQ$ be the path algebra of $Q$ over $\mathbb{F}_q$. Let $P=P(e)$ be the projective module cover the simple module $S_e$. Set $\mathfrak{p}=\underline{\dim} P(e)$. Clearly $\langle\mathfrak{p},\mathfrak{p}\rangle=1=\langle\mathfrak{p},\delta\rangle$ and there exists unique indecomposable preprojective module $L$ with $\underline{\dim}L=\mathfrak{p}+\delta$. Moreover, we have $\textrm{Hom}_{\Lambda}(L,P)=0$ and $\textrm{Ext}_{\Lambda}(L,P)=0$. Let $\mathfrak{C}(P,L)$ be the smallest full subcategory of mod-$\Lambda$ which contains $P$ and $L$ and is closed under taking extensions, kernels of epimorphisms and cokernels of monomorphisms. $\mathfrak{C}(P,L)$ is equivalent to the module category of the Kronecker quiver over $\mathbb{F}_q$. Thus we have an exact embedding $F:\textrm{mod-}K\hookrightarrow \textrm{mod-}\Lambda$, where $K$ is the path algebra of the Kronecker quiver over $\mathbb{F}_q$. We know that the embedding $F$ is independent of the choice of $q$. Hence, this gives rise to an injective homomorphism of algebras $F:\mathcal{H}^{\ast}(K)\rightarrow\mathcal{H}^{\ast}(Q)$. In $\mathcal{H}^{\ast}(K)$, we have defined the element $E_{m\delta_{K}}$ for $m\geq1$. Define $E_{m\delta}=F(E_{m\delta_{K}})$. Since $E_{m\delta_{K}}\in \mathcal{C}^{\ast}(K)$
 , we have $E_{m\delta}\in\mathcal{C}^{\ast}(Q)$.

We may list all non-homogeneous tubes $\mathcal{T}_1,\mathcal{T}_2,\ldots,\mathcal{T}_s$ in mod-$\Lambda$ (in fact $s\leq 3$). For each $\mathcal{T}_i$, let $r_i$ be the period of $\mathcal{T}_i$. For each $\mathcal{T}_i$, as we did in Section 4.1.2, we have the generic composition algebra $\mathcal{C}^{\ast}(\mathcal{T}_i)$ and its integral form $\mathcal{C}^{\ast}(\mathcal{T}_i)_{\mathcal{Z}}$. For each $\mathcal{T}_i$, we have the set $\Pi_{i}^a$ of aperiodic $r_i$-tuples of partitions. We have constructed in Section 4.1.2 the element $E_{\pi_i}$. Then $\{E_{\pi_i}|\pi\in\Pi_{i}^a\}$ is a $\mathcal{Z}$-basis of $\mathcal{C}^{\ast}(\mathcal{T}_i)_{\mathcal{Z}}$.

Let $\mathcal{M}$ be the set of quadruples $\mathbf{c}=(\mathbf{a}_{\mathbf{c}},\mathbf{b}_{\mathbf{c}},\pi_{\mathbf{c}},w_{\mathbf{c}})$ such that $\mathbf{a}_{\mathbf{c}}\in\mathbb{N}^{Prep}_{f}$, $\mathbf{b}_{\mathbf{c}}\in\mathbb{N}^{Prei}_{f}$, $\pi_{\mathbf{c}}=(\pi_{1\mathbf{c}},\pi_{2\mathbf{c}},\ldots,\pi_{s\mathbf{c}})\in\Pi_{1}^a\times\Pi_{2}^a\times\cdots\times\Pi_{s}^a$ and $w_{\mathbf{c}}=({w_1}\geq{w_2}\geq\cdots\geq{w_t})$ is a partition of $m\in\mathbb{Z}_{>0}$.

Then for each $\mathbf{c}\in\mathcal{M}$ we define
\begin{displaymath}
E^{\mathbf{c}}=\langle M(\mathbf{a}_{\mathbf{c}})\rangle\ast E_{\pi_{1\mathbf{c}}}\ast E_{\pi_{2\mathbf{c}}}\ast\cdots\ast E_{\pi_{s\mathbf{c}}}\ast E_{w_{\mathbf{c}}\delta}\ast\langle M(\mathbf{b}_{\mathbf{c}})\rangle,
\end{displaymath}
where $\langle M(\mathbf{a}_{\mathbf{c}})\rangle$ and $\langle M(\mathbf{b}_{\mathbf{c}})\rangle$ are defined in Section 4.1.3, $E_{\pi_{i\mathbf{c}}}$ is defined in Section 4.1.2 and $E_{w_{\mathbf{c}}\delta}$ is defined in Section 4.1.1. Obviously, $\{E^{\mathbf{c}}|\mathbf{c}\in\mathcal{M}\}$ belongs to $\mathcal{C}^{\ast}(Q)$.


\th {Proposition 4.4. [9]}\ {\it The set $\{E^{\mathbf{c}}|\mathbf{c}\in\mathcal{M}\}$ is a $\mathcal{A}$-basis of $\mathcal{C}^{\ast}(Q)_{\mathcal{A}}$.}

From this basis we can get a bar-invariant basis. But it is not the one considered by Lusztig. Hence in [9], another PBW-type basis is constructed. Let us recall it definition.

There is an inner pair $(,)$ on $\mathcal{H}_q^{\ast}(\Lambda)$ defined in [4]. This inner product is also well-defined on $\mathcal{C}^{\ast}(Q)$ which coincides with the pairing defined by Lusztig in [13]. Consider the $\mathbb{Q}(v)$-basis $\{E^{\mathbf{c}}|\mathbf{c}\in\mathcal{M}\}$.
Let $R(\mathcal{C}^{\ast}(Q))$ be the $\mathbb{Q}(v)$-subspace of $\mathcal{C}^{\ast}(Q)$ with the basis $\{E_{\pi_{1\mathbf{c}}}\ast E_{\pi_{2\mathbf{c}}}\ast\cdots\ast E_{\pi_{s\mathbf{c}}}\ast E_{w_{\mathbf{c}}\delta}\}$
where $\pi_{\mathbf{c}}=(\pi_{1\mathbf{c}},\pi_{2\mathbf{c}},\ldots,\pi_{s\mathbf{c}})\in\Pi_{1}^a\times\Pi_{2}^a\times\cdots\times\Pi_{s}^a$, and $w_{\mathbf{c}}=({w_1}\geq{w_2}\geq\cdots\geq{w_t})$ is a partition. It is a subalgebra of $\mathcal{C}^{\ast}(Q)$.

Let $R^a(\mathcal{C}^{\ast}(Q))$ be the subalgebra of $R(\mathcal{C}^{\ast}(Q))$ with the basis $\{E_{\pi_{1\mathbf{c}}}\ast E_{\pi_{2\mathbf{c}}}\ast\cdots\ast E_{\pi_{s\mathbf{c}}}|\pi_{\mathbf{c}}=(\pi_{1\mathbf{c}},\pi_{2\mathbf{c}},\ldots,\pi_{s\mathbf{c}})\in\Pi_{1}^a\times\Pi_{2}^a\times\cdots\times\Pi_{s}^a\}$. For $\alpha,\beta\in\mathbb{N}[I]$, we denote $\alpha\leq\beta$ if $\beta-\alpha\in\mathbb{N}[I]$. It follow that $R(\mathcal{C}^{\ast}(Q))_{\beta}=R^a(\mathcal{C}^{\ast}(Q))_{\beta}$ if $\beta<\delta$. Define $\mathcal{F}_{\delta}=\{x|(x,R^a(\mathcal{C}^{\ast}(Q))_{\delta})=0\}$.

In [9], it is prove that
\begin{displaymath}
R(\mathcal{C}^{\ast}(Q))_{\delta}=R^a(\mathcal{C}^{\ast}(Q))_{\delta}\oplus\mathcal{F}_{\delta}
\end{displaymath}
and $\dim\mathcal{F}_{\delta}=1$. By the method of Schmidt orthogonalization, we may set
\begin{displaymath}
E'_{\delta}=E_{\delta}-\sum_{M(\pi_{i\mathbf{c}}),\underline{\dim} M(\pi_{i\mathbf{c}})=\delta, 1\leq i\leq s}a_{\pi_{i\mathbf{c}}}E_{\pi_{i\mathbf{c}}}
\end{displaymath}
satisfying $\mathcal{F}_{\delta}=\mathbb{Q}(v)E'_{\delta}$.

Now let $R(\mathcal{C}^{\ast}(Q))(1)$ be the subalgebra of $R(\mathcal{C}^{\ast}(Q))$ generated by $R^a(\mathcal{C}^{\ast}(Q))$ and $\mathcal{F}_{\delta}$. We have $R(\mathcal{C}^{\ast}(Q))(1)_{\beta}=R(\mathcal{C}^{\ast}(Q))_{\beta}$ if $\beta<2\delta$. Define
\begin{displaymath}
\mathcal{F}_{2\delta}=\{x|(x,R(\mathcal{C}^{\ast}(Q))(1)_{2\delta})=0\}.
\end{displaymath}
Then $\dim\mathcal{F}_{2\delta}=1$ and $R(\mathcal{C}^{\ast}(Q))_{2\delta}=R(\mathcal{C}^{\ast}(Q))(1)_{2\delta}\oplus\mathcal{F}_{2\delta}$.

In general, define
\begin{displaymath}
\mathcal{F}_{n\delta}=\{x|(x,R(\mathcal{C}^{\ast}(Q))(n-1)_{n\delta})=0\}.
\end{displaymath}
Let $R(\mathcal{C}^{\ast}(Q))(n)$ be the subalgebra of $R(\mathcal{C}^{\ast}(Q))$ generated by $R(\mathcal{C}^{\ast}(Q))(n-1)$ and $\mathcal{F}_{n\delta}$. Then $\dim\mathcal{F}_{n\delta}=1$ and $R(\mathcal{C}^{\ast}(Q))_{n\delta}=R(\mathcal{C}^{\ast}(Q))(n-1)_{n\delta}\oplus\mathcal{F}_{n\delta}$. Similarly, choose $E'_{n\delta}$ such that $E_{n\delta}-E'_{n\delta}\in R(\mathcal{C}^{\ast}(Q))(n-1)_{n\delta}$ and $\mathcal{F}_{n\delta}=\mathbb{Q}(v)E'_{n\delta}$ for all $n>0$.

Let $P_{n\delta}=nE'_{n\delta}$. For a partition $w_{\mathbf{c}}=(1^{r_1}2^{r_2}\cdots t^{r_t})$ of $m\in\mathbb{Z}_{>0}$, let $P_{w_{\mathbf{c}}\delta}=P_{1\delta}^{\ast r_1}\ast \cdots \ast P_{t\delta}^{\ast r_t}$. Let $S_{w_{\mathbf{c}}\delta}$ be the Schur functions corresponding to $P_{w_{\mathbf{c}}\delta}$ and
\begin{displaymath}
F^{\mathbf{c}}=\langle M(\mathbf{a}_{\mathbf{c}})\rangle\ast E_{\pi_{1\mathbf{c}}}\ast E_{\pi_{2\mathbf{c}}}\ast\cdots\ast E_{\pi_{s\mathbf{c}}}\ast S_{w_{\mathbf{c}}\delta}\ast\langle M(\mathbf{b}_{\mathbf{c}})\rangle
\end{displaymath}
for $\mathbf{c}\in\mathcal{M}$.

\th {Proposition 4.5. [9]}\ {\it The set $\{F^{\mathbf{c}}|\mathbf{c}\in\mathcal{M}\}$ is an almost orthonormal basis of $\mathcal{C}^{\ast}(Q)\simeq\mathbf{f}$.}

\subsec{4.2\quad PBW-type basis of $\dot{\mathbf{U}}\mathbf{1}_{\lambda}$}

Let $Q$ be a connected tame quiver without oriented cycles. Consider the root category $\mathcal{R}(Q)$ over some finite field $k$. Let $\Lambda=kQ$.
Remember that $\tilde{\mathcal{P}}$ is the set of isomorphism classes of the objects in $\mathcal{R}(Q)$ and ind$(\tilde{\mathcal{P}})$ is the set of isomorphism classes of the indecomposable objects in $\mathcal{R}(Q)$. The set ind$(\tilde{\mathcal{P}})$ can be divided as follow
\begin{displaymath}
\textrm{ind$(\tilde{\mathcal{P}})$}=\mathbb{P}\,\,\dot{\cup}\,\,\mathbb{T}\,\,\dot{\cup}\,\,T(\mathbb{P})\,\,\dot{\cup}\,\,T(\mathbb{T}).
\end{displaymath}
Fix an embedding of mod-$kQ$ into the root category $\mathcal{R}(Q)$, then $\mathbb{P}=Prep(Q)\dot{\cup}T(Prei(Q))$ and $\mathbb{T}$ is the set of all indecomposable regular $Q$ representations consisting of homogeneous tubes and non-homogeneous tubes $\mathcal{T}_1,\mathcal{T}_2,\cdots,\mathcal{T}_{s}$ appearing in mod-$kQ$.

Let $\tilde{\mathcal{M}}$ be the set of
\begin{displaymath}
\tilde{\mathbf{c}}=
(\mathbf{d}_{\tilde{\mathbf{c}}},\pi_{\tilde{\mathbf{c}}},w_{\tilde{\mathbf{c}}},\mathbf{d}'_{\tilde{\mathbf{c}}},\pi'_{\tilde{\mathbf{c}}},w'_{\tilde{\mathbf{c}}})
\end{displaymath}
where
\begin{displaymath}
\mathbf{d}_{\tilde{\mathbf{c}}}\in \mathbb{N}^{\mathbb{P}}_{f},\quad \mathbf{d}'_{\tilde{\mathbf{c}}}\in \mathbb{N}^{T(\mathbb{P})}_{f},
\end{displaymath}
\begin{eqnarray*}
\pi_{\tilde{\mathbf{c}}}=(\pi_{1\tilde{\mathbf{c}}},\pi_{2\tilde{\mathbf{c}}},\ldots,\pi_{s\tilde{\mathbf{c}}})
\in\Pi_{1}^a\times\Pi_{2}^a\times\cdots\times\Pi_{s}^a,\\
\pi'_{\tilde{\mathbf{c}}}=(\pi'_{1\tilde{\mathbf{c}}},\pi'_{2\tilde{\mathbf{c}}},\ldots,\pi'_{s\tilde{\mathbf{c}}})
\in\Pi_{1}^a\times\Pi_{2}^a\times\cdots\times\Pi_{s}^a,
\end{eqnarray*}
and
\begin{eqnarray*}
w_{\tilde{\mathbf{c}}}=({w_1}\geq{w_2}\geq\cdots\geq{w_t}),\\
w'_{\tilde{\mathbf{c}}}=({w'_1}\geq{w'_2}\geq\cdots\geq{w'_{t'}}),
\end{eqnarray*}
are partitions of $m$ and $m'\in\mathbb{Z}_{>0}$ respectively.
$\mathbb{N}^{\mathbb{P}}_{f}$ is the set of all support-finite function $\mathbf{d}:\mathbb{P}\rightarrow\mathbb{N}$ and $\mathbb{N}^{T(\mathbb{P})}_{f}$ is the set of all support-finite function $\mathbf{d}:T(\mathbb{P})\rightarrow\mathbb{N}$.
Note that $\pi_{\tilde{\mathbf{c}}}=(\pi_{1\tilde{\mathbf{c}}},\pi_{2\tilde{\mathbf{c}}},\ldots,\pi_{s\tilde{\mathbf{c}}})$ and $w_{\tilde{\mathbf{c}}}=({w_1}\geq{w_2}\geq\cdots\geq{w_t})$ defined in Section 4.1, come from objects appearing in $\mathbb{T}$, while $\pi'_{\tilde{\mathbf{c}}}=(\pi'_{1\tilde{\mathbf{c}}},\pi'_{2\tilde{\mathbf{c}}},\ldots,\pi'_{s\tilde{\mathbf{c}}})$ and $w'_{\tilde{\mathbf{c}}}=({w'_1}\geq{w'_2}\geq\cdots\geq{w'_{t'}})$ defined in Section 4.1, come from objects appearing in $T(\mathbb{T})$.

Note that the category $\mathcal{R}(Q)$, so the set $\tilde{\mathcal{M}}$, depend only on the underlying graph of Q.
If $Q'$ is another quiver such that $D^b(kQ)\simeq D^b(kQ')$, they give the same $\tilde{\mathcal{M}}$.

Given any symmetric generalized Cartan matrix $A=(a_{ij})_{n\times n}$ of affine type, we consider a quiver $Q$, the quantum enveloping algebra $\mathbf{U}$ and the modified enveloping algebra $\dot{\mathbf{U}}$ corresponding to the Cartan matrix $A=(a_{ij})_{n\times n}$.

Remember that mod-$kQ$ can be embedding into $\mathcal{R}(Q)$ as a full subcategory. Then ind$(\tilde{\mathcal{P}})=\textrm{ind($\mathcal{P}$)}\dot{\cup}\textrm{ind$(T(\mathcal{P}))$}$.

For $\tilde{\mathbf{c}}=
(\mathbf{d}_{\tilde{\mathbf{c}}},\pi_{\tilde{\mathbf{c}}},w_{\tilde{\mathbf{c}}},\mathbf{d}'_{\tilde{\mathbf{c}}},\pi'_{\tilde{\mathbf{c}}},w'_{\tilde{\mathbf{c}}})
\in\tilde{\mathcal{M}}$,
let $\mathbf{d}_1=\mathbf{d}_{\tilde{\mathbf{c}}}|_{\textrm{ind$(T(Prei(Q)))$}}$ and $\mathbf{d}_2=\mathbf{d}_{\tilde{\mathbf{c}}}|_{\textrm{ind$(Prep(Q))$}}$, and we can denote by $\mathbf{d}_{\tilde{\mathbf{c}}}=(\mathbf{d}_1,\mathbf{d}_2)$. Also, let $\mathbf{d}'_1=\mathbf{d}'_{\tilde{\mathbf{c}}}|_{\textrm{ind$(Prei(Q))$}}$ and $\mathbf{d}'_2=\mathbf{d}'_{\tilde{\mathbf{c}}}|_{\textrm{ind$(T(Prep(Q)))$}}$, and we can denote by $\mathbf{d}'_{\tilde{\mathbf{c}}}=(\mathbf{d}'_1,\mathbf{d}'_2)$. Then $\mathbf{c}_1=(\mathbf{d}_2, \pi_{\tilde{\mathbf{c}}},w_{\tilde{\mathbf{c}}},\mathbf{d}'_1)$ and $\mathbf{c}_2=(\mathbf{d}'_2, \pi'_{\tilde{\mathbf{c}}},w'_{\tilde{\mathbf{c}}},\mathbf{d}_1)$ can be regarded as elements in $\mathcal{M}$ and we can denoted by $\tilde{\mathbf{c}}=(\mathbf{c}_1,\mathbf{c}_2)$.

We identify $\mathcal{C}^{\ast}(Q)$ with $\mathbf{f}$ by the correspondence between $u_{[S_i]}$ and $\theta_i$. So the set
\begin{displaymath}
\{F^{\tilde{\mathbf{c}}}_{\lambda}=\langle{M(\mathbf{d}_1)}\rangle^-\langle M(\mathbf{d}_2)\rangle^+ E_{\pi_{1\tilde{\mathbf{c}}}}^+ E_{\pi_{2\tilde{\mathbf{c}}}}^+\cdots E_{\pi_{s\tilde{\mathbf{c}}}}^+ S_{w_{\tilde{\mathbf{c}}}\delta}^+\langle M(\mathbf{d}'_1)\rangle^+\langle M(\mathbf{d}'_2)\rangle^- E_{\pi'_{1\tilde{\mathbf{c}}}}^-E_{\pi'_{2\tilde{\mathbf{c}}}}^-\cdots E_{\pi'_{s\tilde{\mathbf{c}}}}^-S_{w'_{\tilde{\mathbf{c}}}\delta}^-\mathbf{1}_{\lambda}\}
\end{displaymath}
can be regarded as elements in $\dot{\mathbf{U}}\mathbf{1}_{\lambda}$.

We can also consider the following set
\begin{displaymath}
\{F'^{\tilde{\mathbf{c}}}_{\lambda}=F^{\mathbf{c}_1+}\cdot F^{\mathbf{c}_2-}\mathbf{1}_{\lambda}|\tilde{\mathbf{c}}=(\mathbf{c}_1,\mathbf{c}_2),\tilde{\mathbf{c}}\in\tilde{\mathcal{M}}\}.
\end{displaymath}

\th {Lemma 4.1. }\ {\it The set $\{F'^{\tilde{\mathbf{c}}}_{\lambda}|\tilde{\mathbf{c}}\in\tilde{\mathcal{M}}\}$ is a basis of $\dot{\mathbf{U}}\mathbf{1}_{\lambda}$.}

\pf{Proof.} By Proposition 4.5, $\{F^{\mathbf{c}}|\mathbf{c}\in\mathcal{M}\}$ is a basis of $\mathbf{f}$. In ([13], 23.2.1), Lusztig points out that $\dot{\mathbf{U}}$ is a free $\mathbf{f}\otimes\mathbf{f}^{\textrm{opp}}$-module with basis $(\mathbf{1}_{\lambda})_{\lambda\in P}$. So the set
\begin{displaymath}
\{F^{\mathbf{c}_1+}\cdot F^{\mathbf{c}_2-}\mathbf{1}_{\lambda}
|\mathbf{c}_1\in\mathcal{M},\mathbf{c}_2\in\mathcal{M}\}
\end{displaymath}
is a PBW-type basis of $\dot{\mathbf{U}}\mathbf{1}_{\lambda}$. By $\tilde{\mathbf{c}}=(\mathbf{c}_1,\mathbf{c}_2)$, we have the proposition.

\qed

We denote  by $B'_{Q}(\dot{\mathbf{U}}\mathbf{1}_{\lambda})$ the basis $\{F'^{\tilde{\mathbf{c}}}_{\lambda}|\tilde{\mathbf{c}}\in\tilde{\mathcal{M}}\}$.


For $\tilde{\mathbf{c}}=(\mathbf{c}_1,\mathbf{c}_2),\tilde{\mathbf{c}}'=(\mathbf{c}'_1,\mathbf{c}'_2)\in\tilde{\mathcal{M}}$, we define $\tilde{\mathbf{c}}<\tilde{\mathbf{c}}'$ if and only if $tr|F^{\mathbf{c}_1}|\leq tr|F^{\mathbf{c}'_1}|$, $tr|F^{\mathbf{c}_2}|\leq tr|F^{\mathbf{c}'_2}|$ when $tr|F^{\mathbf{c}_1}|\neq tr|F^{\mathbf{c}'_1}|$ or $tr|F^{\mathbf{c}_2}|\neq tr|F^{\mathbf{c}'_2}|$. If $tr|F^{\mathbf{c}_1}|=tr|F^{\mathbf{c}'_1}|$ and $tr|F^{\mathbf{c}_2}|=tr|F^{\mathbf{c}'_2}|$, we define $\tilde{\mathbf{c}}<\tilde{\mathbf{c}}'$ if and only if $\mathbf{c}_1\preceq\mathbf{c}'_1$ and $\mathbf{c}_2\preceq\mathbf{c}'_2$ but $\tilde{\mathbf{c}}\neq\tilde{\mathbf{c}}'$, where $\prec$  is the order on the set $\mathcal{M}$ in [9].

\th {Lemma 4.2. }\ {\it The transition matrix from $B_{Q}(\dot{\mathbf{U}}\mathbf{1}_{\lambda})$ to $B'_{Q}(\dot{\mathbf{U}}\mathbf{1}_{\lambda})$ under the order $<$ defined above is an invertible lower triangular matrix with diagonal entries are powers of $v$ and off-diagonal entries in $\mathcal{A}$.}

\pf{Proof.} For $x,y\in\mathbf{f}$ homogeneous, write
\begin{displaymath}
(r\otimes1)r(x)=\sum x_1\otimes x_2\otimes x_3
\end{displaymath}
with $x_k\in\mathbf{f}$ homogeneous and
\begin{displaymath}
(\bar{r}\otimes1)\bar{r}(y)=\sum y_1\otimes y_2\otimes y_3
\end{displaymath}
with $y_k\in\mathbf{f}$ homogeneous, where $r:\mathbf{f}\rightarrow\mathbf{f}\otimes\mathbf{f}$ is defined by $r(\theta_i)=\theta_i\otimes 1+1\otimes\theta_i$ and $\bar{r}(x)=\overline{r(\bar{x})}$. By Proposition 3.1.7 in [13], the following equality holds in $\mathbf{U}$:
\begin{displaymath}
x^-y^+=\sum(-1)^{tr|x_1|-tr|x_3|}v^{-tr|x_1|+tr|x_3|}(x_1,y_1)K_{|x_1|}y_2^+x_2^-\{x_3,y_3\}K_{-|x_3|},
\end{displaymath}
where $\{x,y\}=\overline{(\bar{x},\bar{y})}$. Since $tr|x_2|\leq tr|x|$ and $tr|x_2|= tr|x|$ if and only if $x_1=x_3=\mathbf{1}$, $tr|y_2|\leq tr|y|$ and $tr|y_2|=tr|y|$ if and only if $y_1=y_3=\mathbf{1}$, we have the following
\begin{displaymath}
x^-y^+\mathbf{1}_{\lambda}\equiv y^+x^-\mathbf{1}_{\lambda}\,\textrm{mod}\,P(tr|x|-1,tr|y|-1).
\end{displaymath}

Let $\tilde{\mathbf{c}}=(\mathbf{c}_1,\mathbf{c}_2)$. Hence,
\begin{eqnarray*}
&&F^{\tilde{\mathbf{c}}}_{\lambda}\\
&=&\langle{M(\mathbf{d}_1)}\rangle^-\langle M(\mathbf{d}_2)\rangle^+ E_{\pi_{1\tilde{\mathbf{c}}}}^+ E_{\pi_{2\tilde{\mathbf{c}}}}^+\cdots E_{\pi_{s\tilde{\mathbf{c}}}}^+ S_{w_{\tilde{\mathbf{c}}}\delta}^+\langle M(\mathbf{b}'_1)\rangle^+\langle M(\mathbf{d}'_2)\rangle^- E_{\pi'_{1\tilde{\mathbf{c}}}}^-E_{\pi'_{2\tilde{\mathbf{c}}}}^-\cdots E_{\pi'_{s\tilde{\mathbf{c}}}}^-S_{w'_{\tilde{\mathbf{c}}}\delta}^-\mathbf{1}_{\lambda}\\
&=&\langle{M(\mathbf{d}_1)}\rangle^-F^{\mathbf{c}_1+}\langle M(\mathbf{d}'_2)\rangle^- E_{\pi'_{1\tilde{\mathbf{c}}}}^-E_{\pi'_{2\tilde{\mathbf{c}}}}^-\cdots E_{\pi'_{s\tilde{\mathbf{c}}}}^-S_{w'_{\tilde{\mathbf{c}}}\delta}^-\mathbf{1}_{\lambda}\\
&\equiv&F^{\mathbf{c}_1+}\langle{M(\mathbf{d}_1)}\rangle^-\langle M(\mathbf{d}'_2)\rangle^- E_{\pi'_{1\tilde{\mathbf{c}}}}^-E_{\pi'_{2\tilde{\mathbf{c}}}}^-\cdots E_{\pi'_{s\tilde{\mathbf{c}}}}^-S_{w'_{\tilde{\mathbf{c}}}\delta}^-\mathbf{1}_{\lambda}\,\,\,\textrm{mod}\,P(m,n)
\end{eqnarray*}
where $$m=tr|F^{\mathbf{c}_1}|$$ and $$n=tr|\langle{M(\mathbf{d}_1)}\rangle\ast\langle M(\mathbf{d}'_2)\rangle\ast E_{\pi'_{1\tilde{\mathbf{c}}}}\ast E_{\pi'_{2\tilde{\mathbf{c}}}}\ast\cdots\ast E_{\pi'_{s\tilde{\mathbf{c}}}}\ast S_{w'_{\tilde{\mathbf{c}}}\delta}|=tr|F^{\mathbf{c}_2}|.$$
Then
$$F^{\tilde{\mathbf{c}}}_{\lambda}
=F^{\mathbf{c}_1+}\langle{M(\mathbf{d}_1)}\rangle^-\langle M(\mathbf{d}'_2)\rangle^- E_{\pi'_{1\tilde{\mathbf{c}}}}^-E_{\pi'_{2\tilde{\mathbf{c}}}}^-\cdots E_{\pi'_{s\tilde{\mathbf{c}}}}^-S_{w'_{\tilde{\mathbf{c}}}\delta}^-\mathbf{1}_{\lambda}+\sum_{\tilde{\mathbf{c}}'<\tilde{\mathbf{c}}}
\tilde{e}_{\tilde{\mathbf{c}}\tilde{\mathbf{c}}'} F'^{\tilde{\mathbf{c}}'}_{\lambda}$$
for some $\tilde{e}_{\tilde{\mathbf{c}}\tilde{\mathbf{c}}'}\in\mathcal{A}$.

From the definition of the order $\prec$ on $\mathcal{M}$,
\begin{eqnarray*}
&&F^{\mathbf{c}_1+}\langle{M(\mathbf{d}_1)}\rangle^-\langle M(\mathbf{d}'_2)\rangle^- E_{\pi'_{1\tilde{\mathbf{c}}}}^-E_{\pi'_{2\tilde{\mathbf{c}}}}^-\cdots E_{\pi'_{s\tilde{\mathbf{c}}}}^-S_{w'_{\tilde{\mathbf{c}}}\delta}^-\mathbf{1}_{\lambda}\\
&=&v^{f}F^{\mathbf{c}_1+}(\langle M(\mathbf{d}'_2)\rangle^- E_{\pi'_{1\tilde{\mathbf{c}}}}^-E_{\pi'_{2\tilde{\mathbf{c}}}}^-\cdots E_{\pi'_{s\tilde{\mathbf{c}}}}^-S_{w'_{\tilde{\mathbf{c}}}\delta}^-\langle{M(\mathbf{d}_1)}\rangle^-+
\sum_{\mathbf{c}''_2\prec\mathbf{c}_2}e_{\mathbf{c}_2\mathbf{c}''_2} F^{\mathbf{c}''_2-})\mathbf{1}_{\lambda}\\
&=&v^{f}F^{\mathbf{c}_1+}(F^{\mathbf{c}_2-}+
\sum_{\mathbf{c}''_2\prec\mathbf{c}_2}e_{\mathbf{c}_2\mathbf{c}''_2} F^{\mathbf{c}''_2-})\mathbf{1}_{\lambda}\\
&=&v^{f}F^{\mathbf{c}_1+}F^{\mathbf{c}_2-}\mathbf{1}_{\lambda}+
\sum_{\mathbf{c}''_2\prec\mathbf{c}_2} v^{f}e_{\mathbf{c}_2\mathbf{c}''_2}F^{\mathbf{c}_1+}F^{\mathbf{c}''_2-}\mathbf{1}_{\lambda}\\
&=&v^{f}F'^{\tilde{\mathbf{c}}}_{\lambda}+
\sum_{\mathbf{c}''_2\prec\mathbf{c}_2} v^{f}e_{\mathbf{c}_2\mathbf{c}''_2}F^{\mathbf{c}_1+}F^{\mathbf{c}''_2-}\mathbf{1}_{\lambda}\\
&=&v^{f}F'^{\tilde{\mathbf{c}}}_{\lambda}+
\sum_{\tilde{\mathbf{c}}''<\tilde{\mathbf{c}},\mathbf{c}''_1=\mathbf{c}_1} \tilde{e}_{\tilde{\mathbf{c}}\tilde{\mathbf{c}}''}F'^{\tilde{\mathbf{c}}''}_{\lambda},
\end{eqnarray*}
where $\tilde{\mathbf{c}}''=(\mathbf{c}''_1,\mathbf{c}''_2)$, $f=(|\langle{M(\mathbf{d}_1)}\rangle|,|\langle M(\mathbf{d}'_2)\rangle\ast E_{\pi'_{1\tilde{\mathbf{c}}}}\ast E_{\pi'_{2\tilde{\mathbf{c}}}}\ast\cdots\ast E_{\pi'_{s\tilde{\mathbf{c}}}}\ast S_{w'_{\tilde{\mathbf{c}}}\delta}|)$ and
$\tilde{e}_{\tilde{\mathbf{c}}\tilde{\mathbf{c}}''}=v^{f}e_{\mathbf{c}_2\mathbf{c}''_2}$.

Hence
$$F^{\tilde{\mathbf{c}}}_{\lambda}
=v^{f}F'^{\tilde{\mathbf{c}}}_{\lambda}+
\sum_{\tilde{\mathbf{c}}'<\tilde{\mathbf{c}}} \tilde{e}_{\tilde{\mathbf{c}}\tilde{\mathbf{c}}'}F'^{\tilde{\mathbf{c}}'}_{\lambda},$$
where $\tilde{e}_{\tilde{\mathbf{c}}\tilde{\mathbf{c}}'}\in \mathcal{A}$.

The proof is finished.

\qed

Then, we have the following proposition

\th {Proposition 4.6.}\ {\it The set $\{F^{\tilde{\mathbf{c}}}_{\lambda}|\tilde{\mathbf{c}}\in\tilde{\mathcal{M}}\}$ is a $\mathcal{A}$-basis of $\dot{\mathbf{U}}\mathbf{1}_{\lambda}$.}

\pf{Proof.} By Lemma 4.1, $\{F'^{\tilde{\mathbf{c}}}_{\lambda}|\tilde{\mathbf{c}}\in\tilde{\mathcal{M}}\}$ is a basis of $\dot{\mathbf{U}}\mathbf{1}_{\lambda}$. Since the transition matrix from $B_{Q}(\dot{\mathbf{U}}\mathbf{1}_{\lambda})$ to $B'_{Q}(\dot{\mathbf{U}}\mathbf{1}_{\lambda})$ under the order $<$ defined above is an invertible lower triangular matrix with diagonal entries are powers of $v$ and off-diagonal entries in $\mathcal{A}$, the set $\{F^{\tilde{\mathbf{c}}}_{\lambda}|\tilde{\mathbf{c}}\in\tilde{\mathcal{M}}\}$ is also a $\mathcal{A}$-basis of $\dot{\mathbf{U}}\mathbf{1}_{\lambda}$.

\qed

We denote  by $B_{Q}(\dot{\mathbf{U}}\mathbf{1}_{\lambda})$ the basis $\{F^{\tilde{\mathbf{c}}}_{\lambda}|\tilde{\mathbf{c}}\in\tilde{\mathcal{M}}\}$.

\subsec{4.3\quad A bar invariant basis of $\dot{\mathbf{U}}\mathbf{1}_{\lambda}$}

Let $Q$ and $\mathcal{R}(Q)$ as before. There is an order $\prec$ on the set $\mathcal{M}$ in [9].
For $\tilde{\mathbf{c}},\tilde{\mathbf{c}}'\in\tilde{\mathcal{M}}$, we define $\tilde{\mathbf{c}}\prec\tilde{\mathbf{c}}'$ if and only if $\mathbf{c}_1\preceq\mathbf{c}'_1$ and $\mathbf{c}_2\preceq\mathbf{c}'_2$ but $\tilde{\mathbf{c}}\neq\tilde{\mathbf{c}}'$, where $\tilde{\mathbf{c}}=(\mathbf{c}_1,\mathbf{c}_2)$ and $\tilde{\mathbf{c}}'=(\mathbf{c}'_1,\mathbf{c}'_2)$.

Recall that for $\mathbf{c}\in\mathcal{M}$, there exist a monomials $\mathbf{m}_{\mathbf{c}}$ on Chevalley generators $u_{[S_i]}$  satisfying
\begin{displaymath}
\mathbf{m}_{\mathbf{c}}=F^{\mathbf{c}}+\sum_{\mathbf{c}'\prec\mathbf{c}}a_{\mathbf{c}\mathbf{c}'}F^{\mathbf{c}'},
\end{displaymath}
where $a_{\mathbf{c}\mathbf{c}'}\in \mathcal{A}$ ([9]).

Let $a=(a_{\mathbf{c}\mathbf{c}'})$ be the transition matrix from $\{F^{\mathbf{c}}|\mathbf{c}\in\mathcal{M}\}$
to $\{\mathbf{m}_{\mathbf{c}}|\mathbf{c}\in\mathcal{M}\}$, where  $a_{\mathbf{c}\mathbf{c}}=1$ and $a_{\mathbf{c}\mathbf{c}'}=0$ unless $\mathbf{c}'\prec\mathbf{c}$. Note that $a$ is unipotent lower triangular matrix.

Let $\bar{a}$ be obtained from $a$ by applying the $\bar{()}$-involution to each elements of $a$.  Since $\overline{\mathbf{m}_{\mathbf{c}}}=\mathbf{m}_{\mathbf{c}}$, we have
\begin{displaymath}
\mathbf{m}_{\mathbf{c}}=\overline{\mathbf{m}_{\mathbf{c}}}=\sum_{\mathbf{c}'}\bar{a}_{\mathbf{c}\mathbf{c}'}\overline{F^{\mathbf{c}'}},
\end{displaymath}
thus
\begin{displaymath}
\overline{F^{\mathbf{c}}}=\sum_{\mathbf{c}'}{\bar{a}^{-1}}_{\mathbf{c}\mathbf{c}'}\mathbf{m}_{\mathbf{c}'}
=\sum_{\mathbf{c}'}\sum_{\mathbf{c}''}{\bar{a}^{-1}}_{\mathbf{c}\mathbf{c}'}{a}_{\mathbf{c}'\mathbf{c}''}F^{\mathbf{c}''}.
\end{displaymath}
Let $h=\bar{a}^{-1}a$, then $h$ is again a unipotent lower triangular matrix, and $\bar{h}=h^{-1}$.
Similarly to the finite case, there exists a unique unipotent lower triangular matrix $d=(d_{\mathbf{c}\mathbf{c}'})$ with off-diagonal entries in $v^{-1}\mathbb{Q}[v^{-1}]$, such that $d=\bar{d}h$. Then the canonical basis of $\mathbf{f}$ is
\begin{displaymath}
\mathcal{E}^{\mathbf{c}}=F^{\mathbf{c}}+\sum_{\mathbf{c}'\prec\mathbf{c}}d_{\mathbf{c}\mathbf{c}'}F^{\mathbf{c}'},
\end{displaymath}
with $d_{\mathbf{c}\mathbf{c}'}\in v^{-1}\mathbb{Q}[v^{-1}]$ ([9]).

Similarly, we can get a bar-invariant basis of $\dot{\mathbf{U}}\mathbf{1}_{\lambda}$ from
\begin{displaymath}
B'_{Q}(\dot{\mathbf{U}}\mathbf{1}_{\lambda})=\{F^{\mathbf{c}_1+}\cdot F^{\mathbf{c}_2-}\mathbf{1}_{\lambda}
|\tilde{\mathbf{c}}\in\tilde{\mathcal{M}},\tilde{\mathbf{c}}=(\mathbf{c}_1,\mathbf{c}_2)\}
\end{displaymath}
and
\begin{displaymath}
\{\mathbf{m}_{\mathbf{c_1}}^+\cdot\mathbf{m}_{\mathbf{c_2}}^-\mathbf{1}_{\lambda}|
\tilde{\mathbf{c}}\in\tilde{\mathcal{M}},\tilde{\mathbf{c}}=(\mathbf{c}_1,\mathbf{c}_2)\}
\end{displaymath}
under the order $\prec$ on $\tilde{\mathcal{M}}$ defined above. We define $\mathbf{m}_{\tilde{\mathbf{c}}\lambda}=\mathbf{m}_{\mathbf{c_1}}^+\cdot\mathbf{m}_{\mathbf{c_2}}^-\mathbf{1}_{\lambda}$ where $\tilde{\mathbf{c}}=(\mathbf{c}_1,\mathbf{c}_2)$.

By the relation
\begin{displaymath}
\mathbf{m}_{\mathbf{c}}=F^{\mathbf{c}}+\sum_{\mathbf{c}'\prec\mathbf{c}}a_{\mathbf{c}\mathbf{c}'}F^{\mathbf{c}'},
\end{displaymath}
we have
\begin{displaymath}
\mathbf{m}_{\mathbf{c}_1}^+=F^{\mathbf{c}_1+}+\sum_{\mathbf{c}'_1\prec\mathbf{c}_1}a_{\mathbf{c}_1\mathbf{c}'_1}F^{\mathbf{c}'_1+}
\end{displaymath}
and
\begin{displaymath}
\mathbf{m}_{\mathbf{c}_2}^-=F^{\mathbf{c}_2-}+\sum_{\mathbf{c}'_2\prec\mathbf{c}_2}a_{\mathbf{c}_2\mathbf{c}'_2}F^{\mathbf{c}'_2-}
\end{displaymath}
in $\mathbf{U}^{\pm}$ respectively.
Hence, we have
\begin{eqnarray*}
\mathbf{m}_{\tilde{\mathbf{c}}\lambda}&=&\mathbf{m}_{\mathbf{c}_1}^+\cdot \mathbf{m}_{\mathbf{c}_2}^-\mathbf{1}_{\lambda} \\
&=&(F^{\mathbf{c}_1}+\sum_{\mathbf{c}'_1\prec\mathbf{c}_1}a_{\mathbf{c}_1\mathbf{c}_1'}F^{\mathbf{c}'_1})^+\cdot
(F^{\mathbf{c}_2}+\sum_{\mathbf{c}'_2\prec\mathbf{c}_2}a_{\mathbf{c}_2\mathbf{c}'_2}F^{\mathbf{c}'_2})^-\mathbf{1}_{\lambda}\\
&=&F^{\mathbf{c}_1+}\cdot F^{\mathbf{c}_2-}\mathbf{1}_{\lambda}+
F^{\mathbf{c}_1+}\cdot\sum_{\mathbf{c}'_2\prec\mathbf{c}_2}a_{\mathbf{c}_2\mathbf{c}'_2}F^{\mathbf{c}'_2-}\mathbf{1}_{\lambda}+\\
&&\sum_{\mathbf{c}'_1\prec\mathbf{c}_1}a_{\mathbf{c}_1\mathbf{c}_1'}F^{\mathbf{c}'_1+}\cdot F^{\mathbf{c}_2-}\mathbf{1}_{\lambda}+
\sum_{\mathbf{c}'_1\prec\mathbf{c}_1}a_{\mathbf{c}_1\mathbf{c}_1'}F^{\mathbf{c}'_1+}
\cdot\sum_{\mathbf{c}'_2\prec\mathbf{c}_2}a_{\mathbf{c}_2\mathbf{c}'_2}F^{\mathbf{c}'_2-}\mathbf{1}_{\lambda}\\
&=&F'^{\tilde{\mathbf{c}}}_{\lambda}
+\sum_{\tilde{\mathbf{c}}'\prec\tilde{\mathbf{c}}}\tilde{a}_{\tilde{\mathbf{c}}\tilde{\mathbf{c}}'}F'^{\tilde{\mathbf{c}}'}_{\lambda},
\end{eqnarray*}
where $\tilde{\mathbf{c}}=(\mathbf{c}_1,\mathbf{c}_2)$, $\tilde{\mathbf{c}}'=(\mathbf{c}'_1,\mathbf{c}'_2)$ and $\tilde{a}_{\tilde{\mathbf{c}}\tilde{\mathbf{c}}'}=a_{\mathbf{c}_1\mathbf{c}_1'}a_{\mathbf{c}_2\mathbf{c}'_2}$.

The same as above, let $\tilde{a}=(\tilde{a}_{\tilde{\mathbf{c}}\tilde{\mathbf{c}}'})$ be the transition matrix from $\{F'^{\tilde{\mathbf{c}}}_{\lambda}|\tilde{\mathbf{c}}\in\tilde{\mathcal{M}}\}$
to $\{\mathbf{m}_{\tilde{\mathbf{c}}\lambda}|\tilde{\mathbf{c}}\in\tilde{\mathcal{M}}\}$, where  $\tilde{a}_{\tilde{\mathbf{c}}\tilde{\mathbf{c}}}=1$ and $a_{\widetilde{\mathbf{c}}\tilde{\mathbf{c}}'}=0$ unless $\tilde{\mathbf{c}}'\prec\tilde{\mathbf{c}}$. Note that $\tilde{a}$ is unipotent lower triangular matrix with off-diagonal entries in $\mathcal{A}$.

Let $\bar{\tilde{a}}$ be obtained from $\tilde{a}$ by applying the $\bar{()}$-involution to each elements of $\tilde{a}$. Since $\overline{\mathbf{m}_{\tilde{\mathbf{c}}\lambda}}=\mathbf{m}_{\tilde{\mathbf{c}}\lambda}$, we have
\begin{displaymath}
\mathbf{m}_{\tilde{\mathbf{c}}\lambda}=\overline{\mathbf{m}_{\tilde{\mathbf{c}}\lambda}}=\sum_{\tilde{\mathbf{c}}'}\bar{a}_{\tilde{\mathbf{c}}\tilde{\mathbf{c}}'}
\overline{F'^{\tilde{\mathbf{c}}'}_{\lambda}},
\end{displaymath}
thus
\begin{displaymath}
\overline{F'^{\tilde{\mathbf{c}}}_{\lambda}}=\sum_{\tilde{\mathbf{c}}'}\bar{\tilde{a}}^{-1}_{\tilde{\mathbf{c}}\tilde{\mathbf{c}}'}\mathbf{m}_{\tilde{\mathbf{c}}'\lambda}
=\sum_{\tilde{\mathbf{c}}'}\sum_{\tilde{\mathbf{c}}''}\bar{\tilde{a}}^{-1}_{\tilde{\mathbf{c}}\tilde{\mathbf{c}}'}{\tilde{a}}_{\tilde{\mathbf{c}}'\tilde{\mathbf{c}}''}
F'^{\tilde{\mathbf{c}}''}_{\lambda}.
\end{displaymath}

Let $\tilde{h}=\bar{\tilde{a}}^{-1}\tilde{a}$, then $\tilde{h}$ is again a unipotent lower triangular matrix, and $\bar{\tilde{h}}=\tilde{h}^{-1}$.
There exists a unique unipotent lower triangular matrix $\tilde{d}=(\tilde{d}_{\tilde{\mathbf{c}}\tilde{\mathbf{c}}'})$ with off-diagonal entries in $v^{-1}\mathbb{Q}[v^{-1}]$
such that $\tilde{d}=\bar{\tilde{d}}\tilde{h}$. Then
we can define a bar-invariant basis of $\dot{\mathbf{U}}\mathbf{1}_{\lambda}$
\begin{displaymath}
\mathcal{E}^{\tilde{\mathbf{c}}}_{\lambda}=
F'^{\tilde{\mathbf{c}}}_{\lambda}+\sum_{\tilde{\mathbf{c}}'\prec\tilde{\mathbf{c}}}\tilde{d}_{\tilde{\mathbf{c}}\tilde{\mathbf{c}}'}F'^{\tilde{\mathbf{c}}'}_{\lambda},
\end{displaymath}
with $\tilde{d}_{\tilde{\mathbf{c}}'\tilde{\mathbf{c}}}\in v^{-1}\mathbb{Q}[v^{-1}]$. We denoted by $\mathcal{B}_{Q}(\dot{\mathbf{U}}\mathbf{1}_{\lambda})$ the above basis.

\th {Theorem 4.1. }\ {\it $\mathcal{B}_{Q}(\dot{\mathbf{U}}\mathbf{1}_{\lambda})=\{\mathcal{E}^{\tilde{\mathbf{c}}}_{\lambda}|\tilde{\mathbf{c}}\in\tilde{\mathcal{M}}\}
        =\{b^+b'^-\mathbf{1}_{\lambda}|b,b'\in \mathbf{B}\}$.}

\pf{Proof.}\ We use the notation in the above subsection.

First, by the definition of $\tilde{a}$, we have
\begin{displaymath}
\tilde{a}_{\tilde{\mathbf{c}}\tilde{\mathbf{c}}'}=a_{\mathbf{c}_1\mathbf{c}_1'}a_{\mathbf{c}_2\mathbf{c}_2'}
\end{displaymath}
where
$\tilde{\mathbf{c}}=(\mathbf{c}_1,\mathbf{c}_2)$, $\tilde{\mathbf{c}}'=(\mathbf{c}'_1,\mathbf{c}'_2)$.
Hence, we have
\begin{displaymath}
\bar{\tilde{a}}_{\tilde{\mathbf{c}}\tilde{\mathbf{c}}'}=\bar{a}_{\mathbf{c}_1\mathbf{c}_1'}\bar{a}_{\mathbf{c}_2\mathbf{c}_2'}.
\end{displaymath}
Note that
\begin{displaymath}
\sum_{\mathbf{c}'}\bar{a}^{-1}_{\mathbf{c}\mathbf{c}'}\bar{a}_{\mathbf{c}'\mathbf{c}}=1.
\end{displaymath}
We have
\begin{eqnarray*}
\sum_{\tilde{\mathbf{c}}'}(\bar{a}^{-1}_{\mathbf{c}_1\mathbf{c}_1'}\bar{a}^{-1}_{\mathbf{c}_2\mathbf{c}_2'})(\bar{\tilde{a}}_{\tilde{\mathbf{c}}'\tilde{\mathbf{c}}})
&=&\sum_{\tilde{\mathbf{c}}'}(\bar{a}^{-1}_{\mathbf{c}_1\mathbf{c}_1'}\bar{a}^{-1}_{\mathbf{c}_2\mathbf{c}_2'})
(\bar{a}_{\mathbf{c}'_1\mathbf{c}_1}\bar{a}_{\mathbf{c}'_2\mathbf{c}_2})\\
&=&\sum_{\mathbf{c}'_1}(\bar{a}^{-1}_{\mathbf{c}_1\mathbf{c}_1'})(\bar{a}_{\mathbf{c}'_1\mathbf{c}_1})
\sum_{\mathbf{c}'_2}(\bar{a}^{-1}_{\mathbf{c}_2\mathbf{c}_2'})(\bar{a}_{\mathbf{c}'_2\mathbf{c}_2})\\
&=&1.
\end{eqnarray*}
Hence, we have
\begin{displaymath}
\bar{\tilde{a}}^{-1}_{\tilde{\mathbf{c}}\tilde{\mathbf{c}}'}=\bar{a}^{-1}_{\mathbf{c}_1\mathbf{c}_1'}\bar{a}^{-1}_{\mathbf{c}_2\mathbf{c}_2'}.
\end{displaymath}
Then
\begin{eqnarray*}
\tilde{h}_{\tilde{\mathbf{c}}\tilde{\mathbf{c}}''}
&=&\sum_{\tilde{\mathbf{c}}'}\bar{\tilde{a}}^{-1}_{\tilde{\mathbf{c}}\tilde{\mathbf{c}}'}\tilde{a}_{\tilde{\mathbf{c}}'\tilde{\mathbf{c}}''}\\
&=&\sum_{\tilde{\mathbf{c}}'}\bar{a}^{-1}_{\mathbf{c}_1\mathbf{c}_1'}\bar{a}^{-1}_{\mathbf{c}_2\mathbf{c}_2'}a_{\mathbf{c}_1'\mathbf{c}_1''}a_{\mathbf{c}_2'\mathbf{c}_2''}\\
&=&\sum_{\mathbf{c}'_1}\bar{a}^{-1}_{\mathbf{c}_1\mathbf{c}_1'}a_{\mathbf{c}_1'\mathbf{c}_1''}
\sum_{\mathbf{c}'_2}\bar{a}^{-1}_{\mathbf{c}_2\mathbf{c}_2'}a_{\mathbf{c}_2'\mathbf{c}_2''}\\
&=&h_{\mathbf{c}_1\mathbf{c}_1''}h_{\mathbf{c}_2\mathbf{c}_2''}.
\end{eqnarray*}

Next, we will check that
\begin{displaymath}
\tilde{d}_{\tilde{\mathbf{c}}\tilde{\mathbf{c}}'}=d_{\mathbf{c}_1\mathbf{c}_1'}d_{\mathbf{c}_2\mathbf{c}_2'}.
\end{displaymath}
By the uniqueness of $\tilde{d}$, we only need to show
\begin{displaymath}
d_{\mathbf{c}_1\mathbf{c}_1''}d_{\mathbf{c}_2\mathbf{c}_2''}
=\sum_{\tilde{\mathbf{c}}'}\overline{d_{\mathbf{c}_1\mathbf{c}_1'}d_{\mathbf{c}_2\mathbf{c}_2'}}\tilde{h}_{\tilde{\mathbf{c}'}\tilde{\mathbf{c}}''}.
\end{displaymath}
We can calculate it directly:
\begin{eqnarray*}
\sum_{\tilde{\mathbf{c}}'}\overline{d_{\mathbf{c}_1\mathbf{c}_1'}d_{\mathbf{c}_2\mathbf{c}_2'}}\tilde{h}_{\tilde{\mathbf{c}'}\tilde{\mathbf{c}}''}
&=&\sum_{\tilde{\mathbf{c}}'}\bar{d}_{\mathbf{c}_1\mathbf{c}_1'}\bar{d}_{\mathbf{c}_2\mathbf{c}_2'}
h_{\mathbf{c}_1'\mathbf{c}_1''}h_{\mathbf{c}_2'\mathbf{c}_2''}\\
&=&\sum_{\mathbf{c}'_1}\bar{d}_{\mathbf{c}_1\mathbf{c}_1'}h_{\mathbf{c}_1'\mathbf{c}_1''}
\sum_{\mathbf{c}'_2}\bar{d}_{\mathbf{c}_2\mathbf{c}_2'}h_{\mathbf{c}_2'\mathbf{c}_2''}\\
&=&d_{\mathbf{c}_1\mathbf{c}_1''}d_{\mathbf{c}_2\mathbf{c}_2''}.
\end{eqnarray*}
Hence, we have
\begin{displaymath}
\tilde{d}_{\tilde{\mathbf{c}}\tilde{\mathbf{c}}'}=d_{\mathbf{c}_1\mathbf{c}_1'}d_{\mathbf{c}_2\mathbf{c}_2'}.
\end{displaymath}

Now, by the definition,
\begin{eqnarray*}
\mathcal{E}^{\tilde{\mathbf{c}}}_{\lambda}&=&F'^{\tilde{\mathbf{c}}}_{\lambda}
+\sum_{\tilde{\mathbf{c}}'\prec\tilde{\mathbf{c}}}
\tilde{d}_{\tilde{\mathbf{c}}\tilde{\mathbf{c}}'}F'^{\tilde{\mathbf{c}}'}_{\lambda}\\
&=&F'^{\tilde{\mathbf{c}}}_{\lambda}
+\sum_{\tilde{\mathbf{c}}'\prec\tilde{\mathbf{c}}}
d_{\mathbf{c}_1\mathbf{c}'_1}d_{\mathbf{c}_2\mathbf{c}'_2}F'^{\tilde{\mathbf{c}}'}_{\lambda}\\
&=&F^{\mathbf{c_1}+}\cdot F^{\mathbf{c_2}-}\mathbf{1}_{\lambda}
+\sum_{\tilde{\mathbf{c}}'\prec\tilde{\mathbf{c}}}
d_{\mathbf{c}_1\mathbf{c}'_1}d_{\mathbf{c}_2\mathbf{c}'_2}F^{\mathbf{c}'_1+}\cdot F^{\mathbf{c}'_2-}\mathbf{1}_{\lambda}\\
&=&(F^{\mathbf{c_1}}
+\sum_{\mathbf{c}'_1\prec\mathbf{c}_1}d_{\mathbf{c}_1\mathbf{c}'_1}F^{\mathbf{c}'_1})^+\cdot
({F^{\mathbf{c_2}}}
+\sum_{\mathbf{c}'_2\prec\mathbf{c}_2}d_{\mathbf{c}_2\mathbf{c}'_2}{F^{\mathbf{c}'_2}})^-\mathbf{1}_{\lambda}\\
&=&{\mathcal{E}^{\mathbf{c}_1}}^+\cdot{\mathcal{E}^{\mathbf{c}_2}}^-\mathbf{1}_{\lambda}.
\end{eqnarray*}

The proof is finished.

\qed

\th {Remark 4.1. }\ {\it Although we use the embedding of mod-$kQ$ into the root category $\mathcal{R}(Q)$ to construct this basis $\mathcal{B}_{Q}(\dot{\mathbf{U}}\mathbf{1}_{\lambda})$, but this theorem show that this basis is independent of the choice of the orientation of $Q$ in fact.
}

\subsec{4.4\quad A parameterization of the canonical basis of $\dot{\mathbf{U}}\mathbf{1}_{\lambda}$}

Let $\dot{\mathbf{U}}=\oplus_{\lambda\in P}\dot{\mathbf{U}}\mathbf{1}_{\lambda}$ be the modified enveloping algebra corresponding to the quiver $Q$ and $\dot{\mathbf{B}}_{\lambda}$ is the canonical basis of $\dot{\mathbf{U}}\mathbf{1}_{\lambda}$.

\th {Theorem 4.2.}\ {\it We have a bijective map
\begin{displaymath}
\Psi_Q:\tilde{\mathcal{M}}\rightarrow{\dot{\mathbf{B}}_{\lambda}}
\end{displaymath}
given by
\begin{displaymath}
\tilde{\mathbf{c}}\mapsto\mathcal{E}^{\mathbf{c}_1}\diamondsuit_{\lambda}\mathcal{E}^{\mathbf{c}_2},
\end{displaymath}
which is the composition of the following two bijection
\begin{eqnarray*}
\tilde{\mathcal{M}}&\rightarrow&\mathcal{B}_{Q}(\dot{\mathbf{U}}\mathbf{1}_{\lambda})\\
\tilde{\mathbf{c}}&\mapsto&\mathcal{E}^{\tilde{\mathbf{c}}}_{\lambda},
\end{eqnarray*}
and
\begin{eqnarray*}
\mathcal{B}_{Q}(\dot{\mathbf{U}}\mathbf{1}_{\lambda})&\rightarrow&\dot{\mathbf{B}}_{\lambda}\\
b^+b'^-\mathbf{1}_{\lambda}&\mapsto&b\diamondsuit_{\lambda}b'.
\end{eqnarray*}
}

\pf{Proof.}\  The first bijection from $\tilde{\mathcal{M}}$ to $\mathcal{B}_{Q}(\dot{\mathbf{U}}\mathbf{1}_{\lambda})$ comes from our construction of $\mathcal{E}^{\tilde{\mathbf{c}}}_{\lambda}$ and the second bijection from $\mathcal{B}_{Q}(\dot{\mathbf{U}}\mathbf{1}_{\lambda})$ to $\dot{\mathbf{B}}_{\lambda}$ comes from Lusztig ([13], Theorem 25.2.1). By Theorem 4.1, $\mathcal{E}^{\tilde{\mathbf{c}}}_{\lambda}
=\mathcal{E}^{\mathbf{c}_1+}\mathcal{E}^{\mathbf{c}_2-}\mathbf{1}_{\lambda}$. So, we have the theorem.

\qed

Note that the set $\tilde{\mathcal{M}}$ depend only on the root category $\mathcal{R}(Q)$, not depend on the embedding of mod-$kQ$ into $\mathcal{R}(Q)$. The all elements in $\tilde{\mathcal{M}}$ give a parameterization of the canonical basis of the modified quantum enveloping algebra by Theorem 4.2.

\newpage

\vspace{3mm}\th{Acknowledgements}
The authors are grateful to the anonymous referees for their helpful
comments and suggestions.

\vskip0.1in \no {\normalsize \bf References}
\vskip0.1in\parskip=0mm \baselineskip 15pt
\renewcommand{\baselinestretch}{1.15}

\footnotesize
\parindent=6mm

 \REF{1\ }Bridgeland T. Quantum groups via Hall algebras of complexes. \emph{Annals of Math.}, to appear

 \REF{2\ }Cramer T. Double Hall algebras and derived equivalences. \emph{Adv. Math.},  {\bf 224}(3): 1097--1120 (2008)

 \REF{3\ }Deng B, Du J, Xiao J. Generic extensions and canonical bases for cyclic quivers. \emph{Canadian J. Math.}, {\bf 59}: 1260--1283 (2007)

 \REF{4\ }Green J A. Hall algebras, hereditary algebras and quantum groups. \emph{Invent. Math.}, {\bf 120}: 361--377 (1995)

 \REF{5\ }Happle D. On the derived category of a finite-dimensional algebra. \emph{Comment. Math. Helv.}, {\bf 62}: 339--389 (1987)

 \REF{6\ }Happle D. Triangulated categories in the representation theory of finite-dimensional algebras. Londou Math. Soc. LNS 119, Cambridge University Press, 1988

 \REF{7\ }Kashiwara M. Crystal bases of modified quantized enveloping algebra. \emph{Duke Math. J.}, {\bf 73}: 383--413 (1994)

 \REF{8\ }Kapranov M. Heisenberg doubles and derived categories. \emph{J. Algebra}, {\bf 202}(2): 712--744 (1998)

 \REF{9\ }Lin Z, Xiao J, Zhang G. Representations of tame quivers and affine canonical bases. \emph{Publ. RIMS Kyoto Univ.}, {\bf 47}: 825--885 (2011)

 \REF{10\ }Lusztig G. Canonical bases arising from quantized enveloping algebra. \emph{J. Amer. Math. Soc.}, {\bf 3}: 447--498 (1990)

 \REF{11\ }Lusztig G. Quivers, perverse sheaves, and the quantized enveloping algebras. \emph{J. Amer. Math. Soc.}, {\bf 4}: 366--421 (1991)

 \REF{12\ }Lusztig G. Canonical bases in tensor product. \emph{Proc. Nat. Acad. Sci. U.S.A.}, {\bf 89}: 8177--8179 (1992)

 \REF{13\ }Lusztig G. Introduction to quantum groups, vol. 110, Progress in Math., Birkhauser, 1993

 \REF{14\ }Peng L, Xiao J. Root categories and simple Lie algebras. \emph{J. Algebra}, {\bf 198}: 19--56 (1997)

 \REF{15\ }Peng L, Xiao J. Triangulated categories and Kac-Moody algebras. \emph{Invent. Math.}, {\bf 140}: 563--603 (2000)

 \REF{16\ }Reineke M. The monoid of families of quiver representations. \emph{Proc. London, Math. Soc}, {\bf 84}: 663--685 (2002)

 \REF{17\ }Ringel C M. The Hall algebra approach to quantum groups. \emph{Aportaciones Matem$\acute{a}$ticas Comunicaciones}, {\bf 15}: 85--114 (1995)

 \REF{18\ }Ringel C M. Hall algebras and quantum groups. \emph{Invent. Math.}, {\bf 101}: 583--592 (1990)

 \REF{19\ }To\"{e}n B. Derived Hall algebras. \emph{Duke Math J}, {\bf 135}(3): 587--615 (2006)

 \REF{20\ }Xiao J. Hall algebra in a root category. Preprint 95-070, Univ. of Bielefeld, 1995

 \REF{21\ }Xiao J, Xu F. Hall algebras associated to triangulated categories. \emph{Duke Math J},  {\bf 143}(2): 357--373 (2008)

 \hml

\begin{center}
\centerline{\psfig{figure=zkxf33.eps}} \centerline{\footnotesize
Fig. 1.\quad }
\end{center}

{\parbox[c]{60mm}\centerline{\psfig{figure=zkxf33.eps}}
\centerline{\footnotesize Fig. 1.\quad }} {}
{\parbox[c]{60mm}
\begin{center}
\footnotesize Table 1\quad \\\vspace{1.5mm} \doublerulesep 0.4pt
\tabcolsep 19pt
\begin{tabular}{\textwidth}{rcccc}
\hline \hline 表的内容 \hline \hline
\end{tabular}
\end{center}}
}

\newpage
\begin{center}
\footnotesize Table 1\quad \\\vspace{1.5mm} \doublerulesep 0.4pt
\tabcolsep 19pt
\begin{tabular*}{\textwidth}{rcccc}
\hline \hline 表的内容 \hline \hline
\end{tabular*}
\end{center}
\zihao{5}

*******************************************做图*********************
\begin{center}
\centerline{\psfig{figure=zkxf33.eps}} \centerline{\footnotesize
Fig. 1.\quad }
\end{center}

\parbox[c]{60mm}{\centerline{\psfig{figure=zkxf33.eps}}
\centerline{\footnotesize Fig. 1.\quad }}
\parbox[c]{60mm}
{ }


\bf {\boldmath

 \begin{enumerate}
 \item[(1)]
\end{enumerate}

\begin{eqnarray}
 \dot {x} = Ax + Bu,  \quad x(0) = x_0,\nonumber\\
 \dot {x} = Ax + Bu,  \quad x(0) = x_0,
\end{eqnarray}

\def\no{\nonumber}\公式居底
\begin{eqnarray}
\end{eqnarray}

\begin{equation}公式齐缝
\end{equation}

\normalsize

\underline